\documentclass[11pt,a4paper, reqno]{amsart}
\usepackage{amsmath,amsfonts,verbatim}
\usepackage{latexsym}
\usepackage{amssymb,leftidx}
\usepackage{extarrows}
\usepackage{overpic}
\usepackage{calrsfs}
\usepackage{xcolor}
\definecolor{purple}{RGB}{128,0,128}

\usepackage{epsfig}
\usepackage{subfigure}
\usepackage{tikz}
\usepackage{bbm}
\usepackage[font=small,labelfont=bf]{caption}

\usepackage[
    colorlinks=true,
    linkcolor=blue,
    citecolor=blue,
    urlcolor=blue
]{hyperref}
\hypersetup{urlcolor=blue, citecolor=blue}
\usepackage{amssymb}
\usepackage{mathrsfs}
\usepackage{amscd}
\usepackage{cite}

\newcommand\R{\mathbb{R}}
\newcommand\C{\mathbb{C}}
\newcommand\Ran{\rm Ran}
\newcommand\Z{\mathbb{Z}}
\newcommand\supp{\rm supp}

\newcommand\p{{\bf p}}

\usepackage{graphicx}
\usepackage{float}

\numberwithin{equation}{section}
\newtheorem{proposition}{Proposition}[section]
\newtheorem{definition}{Definition}[section]
\newtheorem{lemma}{Lemma}[section]
\newtheorem{theorem}{Theorem}[section]

\newtheorem{remark}{Remark}[section]

\begin{document}
\title[Dirac in Aharonov--Bohm magnetic fields]{Dispersive estimates for Dirac equations\\ in Aharonov-Bohm magnetic fields: massive case }

\author{Federico Cacciafesta}
\address{Federico Cacciafesta:
Dipartimento di Matematica, Universit\'a degli studi di Padova, Via Trieste, 63, 35131 Padova PD, Italy}
\email{federico.cacciafesta@unipd.it}

\author{Piero D'Ancona}
\address{Piero D'Ancona:
Department of Mathematics "Guido Castelnuovo", Sapienza University of Rome, Piazzale Aldo Moro 5, 00185 Rome, Italy}
\email{dancona@mat.uniroma1.it}

\author{Zhiqing Yin}
\address{Zhiqing Yin:
Department of Mathematics, Beijing Institute of Technology, Beijing 100081}
\email{zhiqing.yin@bit.edu.cn}
\author{Junyong Zhang}
\address{Junyong Zhang:
Department of Mathematics, Beijing Institute of Technology, Beijing 100081}
\email{zhang\_junyong@bit.edu.cn}

\thanks{%
The authors are partially supported by the MIUR PRIN project 2020XB3EFL, ``Hamiltonian and Dispersive PDEs'', and by the Gruppo Nazionale per l'Analisi Matematica, la Probabilit\`{a} e le loro Applicazioni (GNAMPA), Project CUP E53C23001670001, by the GNAMPA project Fenomeni lineari e non lineari per equazioni di Dirac (2026) and by National Natural Science Foundation of China (12171031) and Beijing Natural Science Foundation(1242011)}%

\keywords{Dispersive estimates; 
Strichartz estimates;
Aharonov-Bohm potential;
Dirac equation}%

\subjclass[2020]{%
42B37,
35Q40,
35Q41
}%

\begin{abstract}
  We study the dispersive properties of the massive Dirac equation perturbed by an Aharonov--Bohm magnetic field, proving in particular dispersive and Strichartz estimates for the associated flow. This paper is the natural continuation of our previous work in the massless case \cite{cacdanyinzha}, and the overall structure of the analysis remains similar. The presence of a nonzero mass term, however, breaks the scaling invariance of the equation and introduces new difficulties at several stages of the argument. In particular, the appearance of the spectral thresholds $\pm m$ leads to a genuinely different low energy regime, which requires a separate analysis of the singular angular component.
\end{abstract}

\maketitle
    \enlargethispage{3pt}




\section{Introduction}

The purpose of this paper is to continue the analysis initiated in \cite{cacdanyinzha} for the two-dimensional Dirac equation in an Aharonov--Bohm magnetic field, by considering here the massive model. Throughout the paper, the mass $m>0$ is fixed. The system we study is
\begin{equation}\label{eq:dirac}
\begin{cases}
\displaystyle
 i\partial_tu+\mathcal{D}_{A}u=0,\quad u(t,x):\mathbb{R}_t\times\mathbb{R}_x^2\rightarrow\mathbb{C}^{2}\\
u(0,x)=u_0(x)
\end{cases}
\end{equation}
where $\mathcal{D}_{A}$ is given by
\begin{equation}\label{op:D}
\begin{split}
\mathcal{D}_{A}&=\begin{pmatrix}
    -m & (i\partial_{1}+A^{1})-i(i \partial_{2}+A^{2})\\
    (i\partial_{1}+A^{1})+i(i \partial_{2}+A^{2}) &  m
  \end{pmatrix}
\end{split}
\end{equation}
and the {\em Aharonov-Bohm potential} 
 $A(x)=(A^1(x),A^2(x))$ is given by
\begin{equation}\label{AB}
A:\R^2\setminus\{(0,0)\}\to\R^2,
\quad
A(x)=\frac{\Phi}{|x|}\left(-\frac{x_2}{|x|},\frac{x_1}{|x|}\right),\, \quad
x=(x_1,x_2)
\end{equation}
(see \cite{AB59}). Since integer shifts of the magnetic flux are related by a gauge transformation, equivalent to a relabelling of the angular modes, the spectral and dispersive properties of the operator depend only on the fractional part of the flux. After this gauge transformation, we may use $\alpha=\Phi-\lfloor\Phi\rfloor$ in place of $\Phi$ in \eqref{AB}
and assume that
\begin{equation*}
    \alpha\in(0,1),
\end{equation*}
the case of integer flux being equivalent to the case
$\alpha=0$.

In the massless case, the Aharonov--Bohm potential has the same homogeneity as the free Dirac operator, so that the equation possesses a natural scaling structure.  The addition of a mass term breaks this invariance and changes the dispersive picture in an essential way.  In particular, as we shall see, the continuous spectrum develops the two thresholds $\pm m$ and, if
\begin{equation*}
        \p=\sqrt{\rho^2-m^2}
\end{equation*}
denotes the radial momentum, the dispersion relation
\begin{equation*}
        \rho=\pm\sqrt{\p^2+m^2}
\end{equation*}
is asymptotically linear in $\p$ at high energies and quadratic near $\p=0$ after removal of the constant oscillation $e^{\pm imt}$.  Thus the high frequency analysis remains closely related to the massless theory, while the threshold region requires a separate argument.

Our aim is to understand how the localized dispersive and Strichartz estimates obtained in \cite{cacdanyinzha} (see also the Erratum \cite{cacdanyinzha-erratum}) are modified by this new low energy behaviour. As in the massless problem, the angular decomposition naturally separates the flow into a regular sector, containing all nonsingular (``{\em regular}") angular modes, and a {\em ``singular"} sector corresponding to the critical mode $k=0$.  The two components
have rather different structures.  On the regular sector, the square of the Dirac operator coincides with the Friedrichs Aharonov--Bohm Hamiltonian plus $m^2$, and the positive and negative spectral branches can therefore be treated through the corresponding Klein--Gordon flow.  On the singular sector, this reduction is no longer available: the generalized eigenfunctions contain negative order Bessel functions and the propagator also contains mixed off diagonal kernels.  We therefore work directly with the spectral representation of the Dirac flow, isolating the threshold contribution from the high frequency regime.  This yields localized time decay estimates, unweighted Strichartz estimates below the integrability threshold imposed by the singular mode, and weighted estimates beyond that range.

Earlier work on dispersive properties of the Dirac equation with an
Aharonov--Bohm field includes local smoothing and weighted Strichartz
estimates \cite{CF}, and generalized Strichartz estimates
\cite{CYZ}. Related results for the scalar Schr\"odinger, wave and
Klein--Gordon equations were obtained in
\cite{FFFP,FZZ,GaoYinZhangZheng}, respectively.
For a complementary analysis of dispersive estimates and threshold
obstructions for two-dimensional Dirac operators with potentials,
see \cite{EG}. We refer to \cite{cacdanyinzha} for further background.

More recently, Cacciafesta, Danesi and S\'er\'e~\cite{CDS} 
obtained massive AB estimates in mixed radial--angular norms
$L_t^pL_{rdr}^qL_\theta^2$,
using a relativistic Hankel transform. 
For $\alpha\ne1/2$, our 
Theorem~\ref{th:strichartz-singular-massive} on the critical
channel enlarges the range of their Theorem~1.2
(e.g., for $\min\{\alpha,1-\alpha\}=1/4$ it includes $(p,q)=(3,7)$, 
while \cite{CDS} requires $p>3$).
Moreover, our results for the regular component give full
spatial norms in a different admissible range, so the results of
the two papers are not directly ordered. We also establish 
localized decay
for the full matrix propagator, including mixed terms, and weighted
singular estimates beyond the unweighted spatial threshold.

\medskip

Before stating our main results, we introduce the required operator domains and notation. As we shall see in Section \ref{sec:pre}, for any given $\alpha$, the operator $\mathcal{D}_A$ initially defined on $D(\mathcal{D}_{A})= \big[C_{c}^{\infty}(\mathbb{R}^{2}\setminus0)\big]^2$ admits a one parameter family of self-adjoint extensions $\mathcal{D}_{A,\gamma}$ indexed by $\gamma\in[0,\pi)$; the choice of the parameter $\gamma$ determines the boundary conditions at the origin of the domain of $\mathcal{D}_{A,\gamma}$. More precisely, given $\gamma\in[0,\pi)$, the elements of $D(\mathcal{D}_{A,\gamma})$ can be represented as

\begin{equation}\label{preliminarydomain}
  f=
   c\sqrt m
    \begin{pmatrix}
      \cos \gamma \cdot K_{\alpha}((\sqrt2 mr) ) \\
      i e^{i\theta}\sin(\gamma+\pi/4)K_{1-\alpha}(\sqrt2 mr)
    \end{pmatrix}
  + \begin{pmatrix}
  \phi_{0}(r) \\
    e^{i\theta}\ \psi_{0}(r)
  \end{pmatrix}
+
  \sum_{k\in \mathbb{Z}\setminus\{0\}}
  \begin{pmatrix}
    e^{ik \theta}\ \phi_{k}(r) \\
    e^{i(k+1)\theta}\ \psi_{k}(r)
  \end{pmatrix},
\end{equation}
where $K_{\alpha}$ is the Macdonald (modified Bessel) function, $c\in \mathbb{C}$ and, for each $k\in\Z$, $\phi_{k},\psi_{k}$ are regular at the origin and such that \begin{equation*}
  \phi_{k},\psi_{k}\in
    \{\phi\in L^{2}(rdr):
    \phi',\phi/r\in L^{2}(rdr),\ 
    \phi\in C([0,+\infty)),\ 
    \phi(0)=0\}
\end{equation*}
(see Theorem~\ref{the:sadjextDAB}); the angular coefficients also satisfy the graph norm summability condition stated there. For $\alpha\ne1/2$, exactly one choice of extension eliminates the more singular component in \eqref{preliminarydomain}. Its remaining singularity is $r^{-\beta}$, where $\beta=\min\{\alpha,1-\alpha\}<1/2$, and its domain is contained in $H^{1/2}(\R^2;\C^2)$. This is the distinguished extension selected by the usual regularity criterion; see \cite{cuesied,dolestlos}. At $\alpha=1/2$, a nonzero $r^{-1/2}$ singularity is not in $H^{1/2}$, so we select a branch explicitly in Definition~\ref{defSA}, without asserting this regularity criterion. Unlike the massless operator, some massive extensions have a bound state in $(-m,m)$. The distinguished extensions considered here have no point spectrum. The details are given in Subsection~\ref{sub:SAext} and Remark~\ref{rem:boundstate}.
\medskip

As we shall see, the estimates (and the strategy of the proofs) will be quite different in the cases $k=0$ and $k\neq 0$ in decomposition \eqref{preliminarydomain}: it is thus useful to introduce the following orthogonal projections on $L^{2}$:
\begin{equation}\label{def-pro}
  P_0:
  L^{2}(\mathbb{R}^{2})^{2}\to   
  L^{2}((0,\infty),rdr)\otimes h_{0}(\mathbb{S}^{1}),
  \qquad
  P_{\bot}=I-P_{0}
\end{equation}
where
\begin{equation*}
h_k(\mathbb S^1)
=\operatorname{span}_{\C}\left\{
\binom{e^{ik\theta}}{0},\binom{0}{e^{i(k+1)\theta}}
\right\}.
\end{equation*}
Thus the two radial spinor coefficients are independent.

We next specify the functional setting. For $s\in\R$, we define the Sobolev spaces
adapted to the distinguished extension by spectral calculus. Here and below,
$D_A$ and $\mathcal D_A$ both abbreviate $\mathcal D_{A,\gamma}$ after the extension
has been fixed, and $H_A^s$ is already a space of $\C^2$-valued functions:
\begin{equation*}
\|f\|_{H_A^s}:=\|\langle \mathcal D_A\rangle^s f\|_{[L^2(\R^2)]^2},
\qquad
\|f\|_{\dot H_A^s}:=\||\mathcal D_A|^s f\|_{[L^2(\R^2)]^2},
\qquad
\langle \mathcal D_A\rangle:=(I+\mathcal D_A^2)^{1/2}.
\end{equation*}
Since the distinguished extension has spectrum in $(-\infty,-m]\cup[m,\infty)$,
these two norms are equivalent for fixed $m>0$, with constants depending on $m$ and $s$.

For a scalar function $\varphi$, we write $\varphi_j(\lambda)=\varphi(2^{-j}\lambda)$.

We state the estimates for the solution group $e^{it\mathcal D_A}$
of~(\ref*{eq:dirac}).
 With this notation in place, we can state our main results.  The first
one concerns the localized time decay of the flow.  As anticipated
above, the regular and singular sectors behave differently: on the
regular component one obtains the usual frequency-localized decay,
whereas the critical channel requires either a spatial weight or a
restriction on the spatial integrability exponent.

 \begin{theorem}[Time decay estimates]\label{thm:disper}  
  Let $\alpha\in(0,1)$, $\varphi\in\mathcal{C}_c^\infty([1,2])$,
  $\tilde\varphi\in\mathcal{C}_c^\infty([1/4,4])$ such that
  $\varphi\tilde \varphi=\varphi$,
  and let $\mathcal{D}_{A,\gamma}$ be the
  self--adjoint extension of $\mathcal{D}_A$ selected in Definition
  \ref{defSA}.
  Then for any $f\in [L^1{(\R^2)}]^2\cap [L^2{(\R^2)}]^2$
  there is a constant $C$, independent of $j\in\Z$, $t\in\R$ and $f$, such that
 \begin{equation}\label{eq:new121}
\begin{aligned}
&\bigl\|
\varphi_j(|\mathcal{D}_{A,\gamma}|)e^{it\mathcal{D}_{A,\gamma}}P_\perp f
\bigr\|_{[L^\infty(\R^2)]^2}\\
&\quad\le C2^{2j}(1+2^j|t|)^{-1/2}
\sum_{\sigma\in\{+,-\}}
\bigl\|
\widetilde\varphi_j(|\mathcal{D}_{A,\gamma}|)\Pi_\perp^\sigma f
\bigr\|_{[L^1(\R^2)]^2},
\end{aligned}
\end{equation}
where $\Pi_\perp^\pm$ are the positive and negative spectral projections
on the regular subspace. For the critical channel,
  \begin{equation}\label{est:dis0}
\left\|
W_j e^{it\mathcal{D}_A}\varphi_j(|\mathcal{D}_A|)P_0W_j f_j
\right\|_{L^\infty}
\lesssim
2^{2j}(1+2^j|t|)^{-1/2}
\|f_j\|_{L^1},
\qquad
f_j=\widetilde\varphi_j(|\mathcal{D}_A|)P_0f,\end{equation}
  where the weight function matrix $W_{j}$ is given by 
\begin{equation}\label{weig-Wj1}
W_j(|x|)=\begin{cases}
\left(\begin{matrix} (1+(2^j|x|)^{-\alpha})^{-1}& 0 \\0 &  1\end{matrix}\quad \right),\quad \alpha\in(0,\frac12];\\
\left(\begin{matrix} 1& 0 \\0 &  (1+(2^j|x|)^{\alpha-1})^{-1}\end{matrix}\quad \right),\quad \alpha\in(\frac12,1).
\end{cases}
\end{equation}

  Moreover, define
  \begin{equation}\label{q-alp}
    \beta:=\min\{\alpha,1-\alpha\},\qquad
    q(\alpha):=\frac2\beta.
  \end{equation}
  Then for $2\leq q<q(\alpha)$ one has the weightless estimate
  \begin{equation}\label{est:disq}
    \begin{split}
    \Big\|e^{it\mathcal{D}_{A,\gamma}}&\varphi(2^{-j}|\mathcal{D}_{A,\gamma}|)P_{0} f(x)\Big\|_{[L^q(\R^2)]^2}\\&
    \leq C2^{2j(1-\frac2q)}(1+2^{j}|t|)^{-\frac12(1-\frac2q)} 
    \|\tilde\varphi(2^{-j}|\mathcal{D}_{A,\gamma}|)P_{0}f\|_{[L^{q'}(\R^2)]^2}.
    \end{split}
  \end{equation}
\end{theorem}


\begin{remark}
The massive Dirac and Klein--Gordon flows have wave behaviour at high energies
and Schr\"odinger behaviour near the thresholds.
Proposition~\ref{prop:regular-dispersive-massive} gives separate estimates in
these two regimes, including the stronger low energy decay $(1+|t|)^{-1}$.
The theorem records a uniform, simpler consequence. All constants may depend
on the fixed mass, flux, exponents and cutoff functions; no uniform limit
$m\to0$ is asserted. The weighted kernel bound also holds with $f_j$ in
\eqref{est:dis0} replaced by any spinor in $L^1\cap L^2$. In particular,
no commutation of $W_j$ with a spectral cutoff is used. Estimates involving
localized or projected $L^1$ norms are understood when those norms are
finite; no $L^1$ boundedness of the spectral projections is assumed.

The spatial ceiling $q(\alpha)$ comes from the critical
singularity $r^{-\beta}$ and the planar radial measure
$r\,dr$. Indeed,
\begin{equation*}
  \int_0^1 |r^{-\beta}|^q r\,dr
  =\int_0^1 r^{1-\beta q}\,dr<\infty
  \quad\Longleftrightarrow\quad q<\frac2\beta.
\end{equation*}
\end{remark}

\medskip

The preceding localized decay estimates are the starting point for the
corresponding spacetime bounds.  We again distinguish the regular and
singular sectors.  For the regular component, the Klein--Gordon
reduction gives a standard family of Strichartz estimates, while the
critical channel retains the spatial restriction dictated by the
singular generalized eigenfunctions and, in the weighted case, the
additional massive threshold constraint described below.

\begin{theorem}[Strichartz estimates for the regular component]
  \label{th:strichartz-regular-massive}
  Let $\alpha\in(0,1)$ and let $\mathcal{D}_{A,\gamma}$ be the
  self--adjoint extension of $\mathcal{D}_A$ selected in Definition
  \ref{defSA}.
  Let $\eta\in[0,1]$, $s\in \mathbb{R}$,
  $(p,q)\in(2,\infty)^{2}$ or $(p,q)=(\infty,2)$
  such that
  \begin{equation*}
    \frac2p+\frac{1+\eta}{q}\le\frac{1+\eta}{2},\qquad 
    s=(2+\eta)\left(\frac12-\frac1q\right)-\frac1p.
  \end{equation*}
  Then there exists a constant $C$ such that, for all
  $f\in H_A^s(\R^2;\C^2)$,
  \begin{equation}\label{eq:strichartz-regular-massive}
    \|e^{it\mathcal D_{A,\gamma}}P_\bot f\|
      _{L_t^p(\R;[L_x^q(\R^2)]^2)}
    \le
    C\|P_\bot f\|_{H_A^s(\R^2;\C^2)}.
  \end{equation}
\end{theorem}

\begin{remark}
  Choosing $\eta=0$ in Theorem
  \ref{th:strichartz-regular-massive} recovers the usual
  wave admissible range, while $\eta=1$ gives the
  Schr\"{o}dinger type range, as is typical of the
  Klein--Gordon equation, but with a suitable regularity loss.
  Compare with \cite{GaoYinZhangZheng} for the corresponding
  scalar phenomenon.
\end{remark}

For the ``singular component", we first state a family of weighted Strichartz estimates.

\begin{theorem}[Weighted Strichartz estimates for the singular component]
\label{prop:weighted-singular-massive}
Let
$\alpha\in(0,1)$ with $\alpha\neq \frac 12$, 
and let $\mathcal{D}_{A,\gamma}$ be the
  self--adjoint extension of $\mathcal{D}_A$ selected in Definition
  \ref{defSA}.
For $\varepsilon>0$,
 define
\begin{equation*}
W_{\varepsilon,\alpha}(r)
=
\begin{cases}
\begin{pmatrix}
(1+r^{-\alpha-\varepsilon})^{-1}&0\\
0&1
\end{pmatrix},
&0<\alpha< \frac12,\\[5mm]
\begin{pmatrix}
1&0\\
0&(1+r^{-(1-\alpha)-\varepsilon})^{-1}
\end{pmatrix},
&\frac12<\alpha<1.
\end{cases}
\end{equation*}
Let $\beta=\min\{\alpha,1-\alpha\}$ and assume
\begin{equation*}
  \frac 2\beta\le q<\infty,
  \qquad
  \frac{2}{1-\beta}<p\le \infty,
  \qquad
  \theta>\frac{\beta-2/q}{\beta+\varepsilon}.
\end{equation*}
Then the following estimate holds
\begin{equation*}
\|W_{\varepsilon,\alpha}^{\theta}e^{it\mathcal{D}_{A,\gamma}}P_0f\|_{L_t^p(\mathbb R;[L_x^q(\mathbb R^2)]^2)}
\le
C_{\varepsilon,\theta}
\|P_0f\|_{H_A^s(\mathbb R^2;\mathbb C^2)},\qquad        s=1-\frac1p-\frac2q.
\end{equation*}

\end{theorem}

\begin{remark}[Massive threshold effect]

\label{rk:threseffect}
The additional condition
\begin{equation*}
        p>r_\alpha=\frac{2}{1-\min\{\alpha,1-\alpha\}}
\end{equation*}
is a genuinely massive threshold effect. Notice that this additional bound is not caused by a new spatial singularity:
the singular Bessel behaviour driven by the generalized eigenfunctions is the same as in the massless problem. What changes is
the time scale near the thresholds $\rho=\pm m$. Let us give some more comments on this.

Let
\begin{equation*}
        \beta:=\min\{\alpha,1-\alpha\},
        \qquad
        \xi:=\sqrt{\rho^2-m^2}.
\end{equation*}
At $\rho=-m$ when $\alpha<1/2$, and at $\rho=m$ when $\alpha>1/2$,
the coefficient of $J_{-\beta}(\xi r)$ does not vanish. On a fixed annulus
$1\le r\le2$, the spatial weight is bounded below and
$J_{-\beta}(\xi r)\sim C_\beta(\xi r)^{-\beta}$.
Since $d|\rho|\simeq_m \xi\,d\xi$, the corresponding threshold model, after writing the
spectral coefficient in $L^2(d\xi)$, is
\begin{equation*}
        \int_0^1 e^{it\sqrt{m^2+\xi^2}}
        \xi^{\frac12-\beta}g(\xi)\,d\xi,
        \qquad g\in L^2(d\xi).
\end{equation*}
Now
\begin{equation*}
        \sqrt{m^2+\xi^2}
        =
        m+\frac{\xi^2}{2m}+O(\xi^4).
\end{equation*}
Thus, if $g$ is supported in $0<\xi\lesssim\delta$, then for
$|t|\lesssim\delta^{-2}$ the phase variation
\begin{equation*}
        t\bigl(\sqrt{m^2+\xi^2}-m\bigr)
\end{equation*}
is bounded on the support, and after removal of the constant phase, its real part has a fixed sign
on a sufficiently small fixed multiple of that time interval for nonnegative data.

Keep the spatial weight and its parameters fixed; the packet width
$\delta>0$ below is independent of the weight parameter $\varepsilon$.
Taking
\begin{equation*}
        g_\delta(\xi)=\delta^{-1/2}\mathbf 1_{(0,\delta)}(\xi),
        \qquad
        \|g_\delta\|_{L^2(d\xi)}\simeq1,
\end{equation*}
one obtains a contribution of size
\begin{equation*}
        \delta^{-1/2}
        \int_0^\delta \xi^{\frac12-\beta}\,d\xi
        \simeq
        \delta^{1-\beta}
\end{equation*}
on a time interval of length $\delta^{-2}$. Hence the $L_t^p$-norm is at least
of size

\begin{equation*}
        \delta^{1-\beta-\frac2p}.
\end{equation*}
Uniform boundedness as $\delta\to0$ therefore forces
\begin{equation*}
        1-\beta-\frac2p\ge0,
        \qquad\text{i.e.}\qquad
        p\ge \frac2{1-\beta}=r_\alpha.
\end{equation*}
This explains the origin of the restriction in the weighted endpoint estimate. The
strict condition $p>r_\alpha$ in the theorem reflects the fact that our proof uses
strong integrability and does not address the endpoint.

By contrast, in the massless case the phase is linear in the spectral variable. A packet
of width $\delta$ has phase variation of order one already on the shorter time
scale $\delta^{-1}$, so the analogous obstruction is only
$p\ge 1/(1-\beta)$, which is automatic in the Strichartz ranges considered there.

Finally, in the unweighted range $q<q(\alpha)=2/\beta$, the condition
\begin{equation*}
        \frac1p+\frac1q<\frac12
\end{equation*}
already implies $p>r_\alpha$. Indeed,
\begin{equation*}
        q<\frac2\beta
        \quad\Longrightarrow\quad
        \frac1q>\frac\beta2,
\end{equation*}
and therefore
\begin{equation*}
        \frac1p<\frac12-\frac1q
        <
        \frac{1-\beta}{2}
        =
        \frac1{r_\alpha}.
\end{equation*}
\end{remark}

We finally consider the singular component without spatial weights.
In this case the negative order Bessel behaviour of the generalized
eigenfunctions imposes an upper restriction on the spatial exponent.
This is the same obstruction already present in the massless problem:
the massive threshold changes the time scale, but not the local
singularity at the origin.  We therefore recover the unweighted
Strichartz estimates in the range $q<q(\alpha)$.

\begin{theorem}[Strichartz estimates for the singular component]
\label{th:strichartz-singular-massive}
Let $\alpha\in(0,1)$ with $\alpha\neq \frac 12$
and let $\mathcal{D}_A$ be the self-adjoint extension selected in Definition \ref{defSA}.
Assume either $(p,q)=(\infty,2)$ or
\begin{equation*}
(p,q)\in(2,\infty]^2,
\qquad
\frac1p+\frac1q<\frac12.
\end{equation*}
Assume in addition that $q<q(\alpha)$, with
$q(\alpha)$ defined in \eqref{q-alp}.
Then there exists a constant $C$ such that, for all $f\in H_A^s(\R^2;\C^2)$,
\begin{equation}\label{eq:strichartz-singular-massive}
\|e^{it\mathcal D_A}P_0f\|_{L_t^p(\R;[L_x^q(\R^2)]^2)}
\le
C\|P_0f\|_{H_A^s(\R^2;\C^2)},\qquad s=1-\frac1p-\frac2q.
\end{equation}
\end{theorem}

Combining Theorems~\ref{th:strichartz-regular-massive}
and~\ref{th:strichartz-singular-massive}, and using $P_\perp+P_0=I$,
we obtain an estimate for the full flow.
Let $\alpha\in(0,1)\setminus\{1/2\}$,
$\beta=\min\{\alpha,1-\alpha\}$, $\eta\in[0,1]$ and
$p,q\in(2,\infty)$ with
\begin{equation*}
  \frac2p+\frac{1+\eta}{q}\le\frac{1+\eta}{2},
  \qquad
  \frac1p+\frac1q<\frac12,
  \qquad
  q<\frac2\beta.
\end{equation*}
Set
\begin{equation*}
  s_\eta=(2+\eta)\left(\frac12-\frac1q\right)-\frac1p,
  \qquad
  s_0=1-\frac1p-\frac2q.
\end{equation*}
Then, for every $f\in H_A^{s_\eta}(\R^2;\C^2)$,
\begin{equation}\label{eq:strichartz-full-flow-massive}
  \|e^{it\mathcal D_A}f\|_{L_t^p(\R;[L_x^q(\R^2)]^2)}
  \le C\|f\|_{H_A^{s_\eta}(\R^2;\C^2)}.
\end{equation}
The energy case $(p,q)=(\infty,2)$ holds as well, with $s_\eta=s_0=0$.
Indeed, $s_\eta-s_0=\eta(1/2-1/q)\ge0$, so the triangle inequality
and the orthogonality of the angular projections in $H_A^{s_\eta}$ give
\begin{equation*}
  \|e^{it\mathcal D_A}f\|_{L_t^pL_x^q}
  \lesssim
  \|P_\perp f\|_{H_A^{s_\eta}}+\|P_0f\|_{H_A^{s_0}} \lesssim 
  \|f\|_{H_A^{s_\eta}}.
\end{equation*}

{\bf Plan of the paper.}
The paper is organized as follows.  In Section~\ref{sec:pre} we recall
the self-adjoint realizations of the massive Dirac Aharonov--Bohm
operator, construct the generalized eigenfunctions and the associated
relativistic Hankel transform, and derive the spectral representation
of the propagator.  Section~\ref{sec:proofdisp} is devoted to the
localized dispersive estimates; the main new point is the analysis of
the spectral thresholds $\rho=\pm m$ in the singular channel.
Finally, Section~\ref{sec:strichartz} establishes the Strichartz estimates,
beginning with the regular component. For the singular component, we
first prove a uniform square function inequality and the unweighted
estimates using localized spacetime bounds. We then obtain the weighted
estimates by interpolation with a weighted $L_t^rL_x^\infty$ estimate.

\section{Self-adjoint extensions, spectral theory and Hankel transform for the Dirac operator}\label{sec:pre}

This section develops the operator domains and spectral representation used
in the proofs. We retain the framework of the massless analysis in
\cite{cacdanyinzha} and explain the changes caused by the mass term. For earlier treatments of self-adjoint realizations and spectral
analysis of Dirac operators with Aharonov--Bohm fields, we also refer
to \cite{gerb,GTV}.

\subsection{Self-adjoint extensions}\label{sub:SAext}

We begin with the partial-wave decomposition needed to characterize the
self-adjoint extensions of $\mathcal D_A$.

By making use of partial wave decomposition, we can decompose the operator $\mathcal{D}_A$ as follows
\begin{equation*}
\mathcal{D}_{A}\simeq
  \bigoplus_{k\in\mathbb Z}d_k
  \qquad
  d_{k}:=
  \begin{pmatrix}
    -m & i\partial_{k+1-\alpha} \\
    i\partial_{\alpha-k} & m
  \end{pmatrix}, 
    \qquad
  I_2:=
  \begin{pmatrix}
    1 & 0 \\
  0 & 1 
  \end{pmatrix}.
\end{equation*}
Equivalently, after expanding spinors
$(\begin{smallmatrix} \phi \\ \psi \end{smallmatrix})
  \in \big[L^{2}(\mathbb{R}^{2})\big]^{2}$
in the form
\begin{equation*}
  \phi=\sum_{k\in \mathbb{Z}}\phi _{k}(r)e^{ik \theta},
  \qquad
  \psi=\sum_{k\in \mathbb{Z}}\psi _{k}(r)e^{i(k+1)\theta},
\end{equation*}
we can write
\begin{equation*}
  \mathcal{D}_{A}
  \sum_{k}
  \begin{pmatrix}
    e^{ik \theta}\ \phi_{k}(r) \\
    e^{i(k+1)\theta}\ \psi_{k}(r)
  \end{pmatrix}
  =
  \sum_{k}
  \begin{pmatrix}
    e^{ik\theta}\bigl(i\partial_{k+1-\alpha}\psi_k(r)-m\phi_k(r)\bigr)\\
    e^{i(k+1)\theta}\bigl(i\partial_{\alpha-k}\phi_k(r)+m\psi_k(r)\bigr)
  \end{pmatrix}
\end{equation*}
where \begin{equation}\label{notation-der}
\partial_{\sigma}=\partial_{r}+\frac{\sigma}{r},
  \qquad
  \partial_{r}=\frac{d}{dr},
  \qquad
  \sigma\in \mathbb{R}.
\end{equation}
For each $k$, the minimal radial operator has domain
$[C_c^\infty((0,\infty))]^2$. Its closure satisfies
\begin{equation*}
  D(\overline{d}_{k})=X^{2}=X \times X:=Y
\end{equation*}
where
\begin{equation}\label{eq:domain}
  X=
  \{\phi\in L^{2}(rdr):
 \phi',\phi/r\in L^{2}(rdr),\ 
  \phi\in C([0,+\infty)),\ 
  \phi(0)=0\}
\end{equation}
We are now in a position to state our result.
\begin{theorem}[Self-adjoint extensions of the Dirac AB operator]\label{the:sadjextDAB}
  Any self-adjoint extension of $\mathcal{D}_{A}$ is of the
  form $\mathcal{D}_{A, \gamma}$, where $\gamma\in[0,\pi)$,
  \begin{equation*}
    D(\mathcal{D}_{A,\gamma})\simeq
    Y^\gamma\oplus\bigoplus_{k\ne0}Y,
  \end{equation*}
where
\begin{equation}\label{Yspace}
Y^\gamma=Y\dotplus\operatorname{span}_{\C}
  \left\{
  \begin{pmatrix}
  \cos\gamma\,K_\alpha(\sqrt2 mr)\\
  i\sin(\gamma+\pi/4)K_{1-\alpha}(\sqrt2 mr)
  \end{pmatrix}
  \right\}.
\end{equation}
The direct sum is understood in the graph norm: writing
$f_k=(\phi_k,\psi_k)^T$, one requires
\begin{equation}\label{eq:graph-summability}
\sum_{k\in\mathbb Z}
\left(\|f_k\|_{L^2(rdr)^2}^2+
\|d_k^\gamma f_k\|_{L^2(rdr)^2}^2\right)<\infty,
\qquad d_k^\gamma:=\overline d_k\quad(k\ne0).
\end{equation}
The operator itself is
  \begin{equation*}
    \mathcal{D}_{A,\gamma}\simeq
d_0^\gamma\oplus\bigoplus_{k\ne0}\overline d_k.
  \end{equation*}
  Any $f\in D(\mathcal{D}_{A,\gamma})$ can be written as
  \begin{equation*}
    f=
    \sum_{k\in \mathbb{Z}}
    \begin{pmatrix}
      e^{ik \theta}\ \phi_{k}(r) \\
      e^{i(k+1)\theta}\ \psi_{k}(r)
    \end{pmatrix}
  \end{equation*}
  with $(\begin{smallmatrix}\phi_{k}\\\psi_{k}\end{smallmatrix})\in Y$
  for $k\neq0$ while
$(\begin{smallmatrix}\phi_{0}\\\psi_{0}\end{smallmatrix})
    \in Y^{\gamma}$
  for $k=0$, i.e., for some $c\in \mathbb{C}$ and
$(\begin{smallmatrix}\tilde\phi_{0}\\\tilde\psi_{0}
    \end{smallmatrix}) \in Y$,
  \begin{equation}\label{domdec}
    \begin{pmatrix}
      \phi_{0} \\
      \psi_{0}
    \end{pmatrix}=
    c\sqrt m
    \begin{pmatrix}
      \cos \gamma \cdot K_{\alpha}((\sqrt2 mr) ) \\
      i\sin (\gamma+\pi/4)\cdot K_{1-\alpha}((\sqrt2 mr) ) 
    \end{pmatrix}
    +\begin{pmatrix}
       \tilde\phi_{0}(r) \\
       \tilde\psi_{0}(r)
    \end{pmatrix}
    =f_{sing}+\begin{pmatrix}
       \tilde\phi_{0}(r) \\
       \tilde\psi_{0}(r)
    \end{pmatrix}.
  \end{equation}
  On the regular part, whose radial coefficients lie in $Y$ and satisfy
\eqref{eq:graph-summability}, the operator acts as
  \begin{equation}\label{DAf}
    \mathcal{D}_{A,\gamma}f=
    \sum_{k\in\Z}
    \begin{pmatrix}
    e^{ik \theta}\ (i\partial_{k+1-\alpha}\psi_{k}(r) -m \phi_k(r))\\
    e^{i(k+1)\theta}\ (i\partial_{\alpha-k}\phi_{k}(r)+m \psi_k(r))
    \end{pmatrix}.
    \end{equation}
    while its radial action on $f_{sing}$ in \eqref{domdec} is
    \begin{equation}\label{action}
d_0^\gamma f_{sing}=cm\sqrt m
\begin{pmatrix}
    \sin \gamma \cdot K_{\alpha}(\sqrt2 mr) \\
    -i\cos(\gamma+\pi/4)K_{1-\alpha}(\sqrt2 mr)
  \end{pmatrix}.
\end{equation}
\end{theorem}

\begin{proof}
The graph norm of the minimal radial operator is determined by
\begin{equation*}
\begin{aligned}
\|d_k f\|_{L^2(rdr)^2}^2
&=\|f'\|_{L^2(rdr)^2}^2+m^2\|f\|_{L^2(rdr)^2}^2\\
&\quad+\|(k-\alpha)f_1/r\|_{L^2(rdr)}^2
+\|(k+1-\alpha)f_2/r\|_{L^2(rdr)}^2.
\end{aligned}
\end{equation*}
for $f\in[C_c^\infty((0,\infty))]^2$. In particular, the mass term must be
included in this identity, and the closure domain is $Y$.
The deficiency spaces at $\pm im$ vanish for $k\ne0$. For $k=0$ they are
one-dimensional, with generators
\begin{equation*}
u_+(r)=\binom{K_\alpha(\sqrt2 mr)}{-e^{-i\pi/4}K_{1-\alpha}(\sqrt2 mr)},
\qquad
u_-(r)=\binom{K_\alpha(\sqrt2 mr)}{e^{i\pi/4}K_{1-\alpha}(\sqrt2 mr)},
\end{equation*}
respectively. Indeed, $d_0^*u_\pm=\pm im u_\pm$, and both components are
in $L^2(rdr)$ precisely in the critical channel. These two generators have
equal norms. Von Neumann's extension theorem therefore gives the additional
domain vector
\begin{equation*}
\tfrac12(e^{-i\gamma}u_++e^{i\gamma}u_-)
=\binom{\cos\gamma\,K_\alpha(\sqrt2 mr)}
{i\sin(\gamma+\pi/4)K_{1-\alpha}(\sqrt2 mr)},
\end{equation*}
with $\gamma$ taken modulo $\pi$. Applying $d_0^*$ gives \eqref{action}.
This also explains the argument $\sqrt2\,mr$, rather than $\sqrt{2mr}$,
in the modified Bessel functions.
\end{proof}

Given Theorem above, we can define a {\em distinguished} self-adjoint extension of $\mathcal{D}_A$ as follows:

\begin{definition}\label{defSA}
For $\alpha\in(0,1)$, we select the following self-adjoint extension from
Theorem~\ref{the:sadjextDAB}, called the distinguished extension in this paper:
\begin{itemize}
\item If $\alpha\in(0,1/2]$: we take $\gamma=\frac{3\pi}4$;
\item if $\alpha\in(1/2,1)$: we take $\gamma=\frac\pi2$.
\end{itemize}
\end{definition} 

For $\alpha\ne1/2$, this is the extension selected by the
$H^{1/2}$ regularity criterion. Indeed, $K_\sigma$ decays exponentially at
infinity, and
\begin{equation*}
K_\sigma(r)\sim 2^{\sigma-1}\Gamma(\sigma)r^{-\sigma}
\quad(r\downarrow0),\qquad 0<\sigma<1.
\end{equation*}
A localized singularity $r^{-\sigma}$ belongs to $H^{1/2}(\R^2)$ exactly
when $\sigma<1/2$. At $\alpha=1/2$, the chosen branch is the limit from
$\alpha<1/2$; it does not satisfy this $H^{1/2}$ criterion. The selected
extensions have purely absolutely continuous spectrum, as follows from
Proposition~\ref{pro:hankel} and the corresponding construction with the
components interchanged.

\subsection{Spectral theory and generalized eigenfunctions}\label{sec:3}
We now briefly review the spectral theory as well as the construction of the generalized eigenfunctions for the Dirac Aharonov--Bohm model in the massive case. There are two main differences from the massless case. The mass breaks
scaling invariance, so the spatial and energy variables no longer enter
the eigenfunctions through a single scaling parameter. It also allows
point spectrum, depending on the self-adjoint extension. More precisely, the critical angular channel may contain a
unique eigenvalue in the spectral gap $(-m,m)$. As we shall see below,
for $\gamma\in[0,\pi)$ such an eigenvalue exists precisely for

\begin{equation*}
    \gamma\in[0,\pi/2)\cup(3\pi/4,\pi),
\end{equation*}
whereas no eigenvalue in the gap is present for
\begin{equation*}
    \gamma\in[\pi/2,3\pi/4].
\end{equation*}
In particular, the distinguished extensions introduced above have no
point spectrum in the gap.

All generalized eigenfunction formulas below are understood for $|\rho|>m$;
values at the thresholds do not affect the spectral integrals.
Let us now explicitly compute the generalized eigenfunctions of the operator; since these will have a key role in the rest of our paper, and since the massive case exhibits some significant, structural differences from the massless case, we provide some fairly complete computations.

The eigenfunctions of the self-adjoint operator $d_{k}^\gamma$ are found as the solutions to the following equation\footnote{ Notice that we cannot, as we did in the massless case, simply solve the eigenvalue equation for $\rho=1$ and then obtain the other solutions by scaling, as the presence of a mass term destroys the scaling of the equation.}:
\begin{equation}\label{eq:geneg}
d_{k}^\gamma \Psi_k(\rho, r)=d_{k}^\gamma \left(\begin{matrix} \psi_k^1(\rho, r)\\ \psi_k^2(\rho, r)\end{matrix}\right)=\rho  \left(\begin{matrix} \psi_k^1(\rho, r)\\ \psi_k^2(\rho, r)\end{matrix}\right)
\end{equation}
where we recall the action of the operator $d_k^\gamma$ is given by
 \begin{equation*}
 d_{k}^{\gamma}=
  \begin{pmatrix}
    -m & i\partial_{k+1-\alpha} \\
    i\partial_{\alpha-k} & m
  \end{pmatrix}=
  \begin{pmatrix}
    -m & i \partial_{r}+i\frac{k+1-\alpha}r \\
    i\partial_{r}+i\frac{\alpha-k} r& m 
  \end{pmatrix}
 \end{equation*}
so that, explicitly,
\begin{equation}\label{explicitsyst}
\begin{cases}
-m\psi^1_k+\left(i\partial_{r}+i\frac{k+1-\alpha}r\right)\psi^2_k=\rho \psi^1_k,
\\ \left(i\partial_{r}+i\frac{\alpha-k}r\right)\psi^1_k+m\psi_k^2=\rho \psi^2_k.
  \end{cases}
\end{equation}
Each component of a solution of \eqref{eq:geneg} satisfies the following
second-order equation. The first-order relations \eqref{explicitsyst} must
still be imposed to link the two components:

\begin{equation}\label{sys-bessel'}
\begin{cases}
  \frac{d^2}{dr^2} \psi^1_k+\frac1r  \frac{d}{dr} \psi^1_k+((\rho^2-m^2)-\frac{(k-\alpha)^2}{r^2})\psi^1_k=0,\\
    \frac{d^2}{dr^2} \psi^2_k+\frac1r  \frac{d}{dr} \psi^2_k+((\rho^2-m^2)-\frac{(k+1-\alpha)^2}{r^2})\psi^2_k=0.
  \end{cases}
\end{equation}
These are Bessel equations.
For $\rho^2>m^2$, the scalar solution bases are $J_\nu$ and $Y_\nu$;
for $\rho^2<m^2$, they are the modified Bessel functions $I_\nu$ and $K_\nu$.
We distinguish these two cases.
\medskip

$\bullet$ {\em Case $\rho^2>m^2$.} Solving the system \eqref{sys-bessel'} for $\rho>0$ gives, with the requirement of square integrability near the origin in the case $k\neq0$,
\begin{equation}\label{gen-eig'}
\begin{aligned}
\psi_k^1(\rho,r)&=\sqrt{\rho-m}\,J_{|k-\alpha|}(\p r),\qquad k\ne0,\\
\psi_k^2(\rho,r)&=
\begin{cases}
-i\sqrt{\rho+m}\,J_{k+1-\alpha}(\p r),&k\ge1,\\
i\sqrt{\rho+m}\,J_{\alpha-k-1}(\p r),&k\le-1.
\end{cases}
\end{aligned}
\end{equation}
where we are denoting with $\p=\sqrt{\rho^2-m^2}$ .

In the remaining case $k=0$, both branches in \eqref{geneige0} are square integrable near the origin. Therefore the choice of the generalized eigenfunction is dictated by the boundary condition which characterizes the self-adjoint extension in Theorem \ref{the:sadjextDAB}. The convenient way to encode such a condition is through the singular coefficients of the two spinor components.

For a spinor $f=(\phi,\psi)^T$ such that the limits below exist, we define
\begin{equation}\label{bdy}
\Gamma_1(f):=\frac{(\sqrt2 m)^\alpha}{2^{\alpha-1}\Gamma(\alpha)}\lim_{r\to0}r^\alpha\phi(r),
\qquad
\Gamma_2(f):=\frac{(\sqrt2 m)^{1-\alpha}}{2^{-\alpha}\Gamma(1-\alpha)}\lim_{r\to0}r^{1-\alpha}\psi(r).
\end{equation}
Using the asymptotics
\begin{equation*}
K_\nu(z)\sim 2^{\nu-1}\Gamma(\nu)z^{-\nu},\qquad z\to0,\quad \nu\in(0,1),
\end{equation*}
and the domain decomposition \eqref{domdec}, one obtains that for every $f\in D(d_0^\gamma)$,
\begin{equation}\label{eq:critical-bc}
\binom{\Gamma_1(f)}{\Gamma_2(f)}
=C\binom{\cos\gamma}{i\sin(\gamma+\pi/4)}
\end{equation}
for some constant $C\in\C$.
Equivalently, there exists $c\in\C$ such that
\begin{equation*}
\Gamma_1(f)=c\sqrt m\cos\gamma,
\qquad
\Gamma_2(f)=ic\sqrt m\sin(\gamma+\pi/4).
\end{equation*}
Accordingly, for $\rho>m$ and $k=0$ we start from the general solution
\begin{equation}\label{geneige0}
\Psi_0(\rho,r)=
\begin{pmatrix}
 a(\rho)\sqrt{\rho-m}\,J_{-\alpha}(\p r)+b(\rho)\sqrt{\rho-m}\,J_{\alpha}(\p r)\\[1mm]
 -i a(\rho)\sqrt{\rho+m}\,J_{1-\alpha}(\p r)+i b(\rho)\sqrt{\rho+m}\,J_{\alpha-1}(\p r)
\end{pmatrix}.
\end{equation}
Using the asymptotics of $J_\nu$ at the origin, namely
\begin{equation*}
J_\nu(z)\sim \frac1{\Gamma(\nu+1)}\Bigl(\frac z2\Bigr)^\nu,
\qquad z\to0,
\end{equation*}
one computes
\begin{equation*}
\Gamma_1\bigl(\Psi_{0}(\rho,\cdot)\bigr)
=\frac{2}{\Gamma(\alpha)\Gamma(1-\alpha)}
\,a(\rho)\sqrt{\rho-m}\Bigl(\frac{\sqrt2 m}{\p}\Bigr)^\alpha,
\end{equation*}
\begin{equation*}
\Gamma_2\bigl(\Psi_{0}(\rho,\cdot)\bigr)
=i\frac{2}{\Gamma(\alpha)\Gamma(1-\alpha)}
\,b(\rho)\sqrt{\rho+m}\Bigl(\frac{\sqrt2 m}{\p}\Bigr)^{1-\alpha}.
\end{equation*}
Hence $\Psi_0(\rho,\cdot)$ satisfies the boundary condition at the origin if and only if these two quantities are proportional to $(\cos\gamma,i\sin(\gamma+\pi/4))$. In particular, whenever $\cos\gamma\neq0$, one must have
\begin{equation}\label{eq:critratio}
\frac{b(\rho)}{a(\rho)}
=
\frac{\sin(\gamma+\pi/4)}{\cos\gamma}
\sqrt{\frac{\rho-m}{\rho+m}}
\Bigl(\frac{\p}{\sqrt2 m}\Bigr)^{1-2\alpha}.
\end{equation}
These generalized eigenfunctions are not globally square integrable; only their
boundary behaviour is being prescribed. If $\cos\gamma=0$, then $a(\rho)=0$; if $\sin(\gamma+\pi/4)=0$, then necessarily $b(\rho)=0$.

\begin{remark}\label{rem:distinguished-critical}
For the distinguished self-adjoint extensions introduced in Definition \ref{defSA}, the critical generalized eigenfunction simplifies drastically.
\begin{itemize}
\item If $\alpha\in(0,1/2]$, then $\gamma=3\pi/4$ and therefore $\sin(\gamma+\pi/4)=0$. Hence $b(\rho)=0$, and we may choose
\begin{equation*}
\Psi_{0}(\rho, r)=
\frac{1}{\sqrt{2}} 
\begin{pmatrix}
\sqrt{\rho-m}J_{-\alpha}(\p r)\\
-i\sqrt{\rho+m}J_{1-\alpha}(\p r)
\end{pmatrix},
\qquad \rho\ge m.
\end{equation*}
\item If $\alpha\in(1/2,1)$, then $\gamma=\pi/2$ and therefore $\cos\gamma=0$. Hence $a(\rho)=0$, and one may take instead the branch built on $(J_{\alpha},J_{\alpha-1})$.
\end{itemize}
\end{remark}

\begin{remark}\label{negenrk}
The negative-energy case is completely analogous. Since in the sequel we only need the distinguished extension corresponding to $\alpha\in(0,1/2]$, we simply record here the formulas for $\rho\le -m$ in that case:
 \begin{equation}\label{neggen-eig'}
\begin{aligned}
\psi_k^1(\rho,r)&=\sqrt{|\rho|+m}\,J_{|k-\alpha|}(\p r),\qquad k\ne0,\\
\psi_k^2(\rho,r)&=
\begin{cases}
i\sqrt{|\rho|-m}\,J_{k+1-\alpha}(\p r),&k\ge1,\\
-i\sqrt{|\rho|-m}\,J_{\alpha-k-1}(\p r),&k\le-1.
\end{cases}
\end{aligned}
\end{equation}
For $\rho<-m$ and $k\neq0$, the normalized generalized eigenfunction is
\begin{equation*}
\Psi_k(\rho,r)
:=
\frac1{\sqrt2}
\begin{pmatrix}
\psi_k^1(\rho,r)\\
\psi_k^2(\rho,r)
\end{pmatrix}.
\end{equation*}
and for $k=0$
 \begin{equation}\label{negeq:psik0}
    \Psi_{0}(\rho, r)=
    \frac{1}{\sqrt{2}} 
    \begin{pmatrix}
  \sqrt{|\rho|+m}  J_{-\alpha}(\p r)\\
  i \sqrt{|\rho|-m}    J_{1-\alpha}(\p r)
    \end{pmatrix}
    \qquad {\rm for}\quad \rho\le -m.
  \end{equation}
For a general extension $\gamma$, the critical channel $k=0$ again involves an energy-dependent linear combination of the two singular branches, exactly as in \eqref{geneige0}--\eqref{eq:critratio}.
\end{remark}
$\bullet$ {\em Case $\rho^2<m^2$.}
Let $-m<\rho<m$ and denote
\begin{equation*}
\kappa=\sqrt{m^2-\rho^2}.
\end{equation*}
Since $K_\nu\in L^2(rdr)$ if and only if $|\nu|<1$, for $\alpha\in(0,1)$ a possible bound state may only occur in the critical channel $k=0$. In that channel the decaying solution is unique up to a multiplicative constant and is given by
\begin{equation}\label{bdstatesol}
\Psi^{E_b}_0(\rho, r)=C_N\begin{pmatrix}
  \sqrt{m-\rho}\,K_{\alpha}(\kappa r)\\
  i\sqrt{m+\rho}\,K_{1-\alpha}(\kappa r)
  \end{pmatrix}.
\end{equation}
Applying the boundary maps \eqref{bdy} yields
\begin{equation*}
\Gamma_1\bigl(\Psi^{E_b}_0(\rho,\cdot)\bigr)
=C_N\sqrt{m-\rho}\Bigl(\frac{\sqrt2 m}{\kappa}\Bigr)^\alpha,
\qquad
\Gamma_2\bigl(\Psi^{E_b}_0(\rho,\cdot)\bigr)
=iC_N\sqrt{m+\rho}\Bigl(\frac{\sqrt2 m}{\kappa}\Bigr)^{1-\alpha}.
\end{equation*}
Therefore $\Psi^{E_b}_0(\rho,\cdot)$ belongs to $D(d_0^\gamma)$ if and only if
\begin{equation}\label{bdstateen}
\binom{\sqrt{m-\rho}\left(\frac{\sqrt2 m}{\kappa}\right)^\alpha}
{i\sqrt{m+\rho}\left(\frac{\sqrt2 m}{\kappa}\right)^{1-\alpha}}
=C
\binom{\cos\gamma}{i\sin(\gamma+\pi/4)}
\end{equation}
for some constant $C\in\C$. 
Equivalently, whenever $\cos\gamma\neq0$,
\begin{equation*}
\frac{\sin(\gamma+\pi/4)}{\cos\gamma}
=
\sqrt{\frac{m+\rho}{m-\rho}}
\Bigl(\frac{\sqrt2 m}{\kappa}\Bigr)^{1-2\alpha}.
\end{equation*}
This gives the implicit equation for the bound-state energy. The
existence and uniqueness of a solution in the gap, as a function of
the extension parameter $\gamma$, are described in Remark~\ref{rem:boundstate}
below. In particular, the distinguished choices
$\gamma=\pi/2$ and $\gamma=3\pi/4$ do not produce any eigenvalue in
$(-m,m)$.

\medskip
\begin{remark}
\label{rem:boundstate}
For a general self-adjoint extension, the critical channel may produce
a unique eigenvalue in the spectral gap $(-m,m)$. Indeed, for every
$E\in(-m,m)$ there exists, up to multiplication by a scalar, a unique
solution which decays at infinity, and its boundary coefficients satisfy
\begin{equation*}
    \frac{\Gamma_2(\Psi_E)}
         {i\Gamma_1(\Psi_E)}
    =
    R_\alpha(E),
\end{equation*}
where
\begin{equation*}
    R_\alpha(E)
    :=
    \sqrt{\frac{m+E}{m-E}}
    \left(
        \frac{\sqrt{2}\,m}{\sqrt{m^2-E^2}}
    \right)^{1-2\alpha}.
\end{equation*}
Equivalently,
\begin{equation*}
    R_\alpha(E)
    =
    (\sqrt{2}\,m)^{1-2\alpha}
    (m+E)^\alpha
    (m-E)^{-(1-\alpha)}.
\end{equation*}
In particular,
\begin{equation*}
    \frac{d}{dE}\log R_\alpha(E)
    =
    \frac{\alpha}{m+E}
    +
    \frac{1-\alpha}{m-E}
    >0,
\end{equation*}
and
\begin{equation*}
    \lim_{E\to-m^+}R_\alpha(E)=0,
    \qquad
    \lim_{E\to m^-}R_\alpha(E)=+\infty.
\end{equation*}
Thus $R_\alpha$ is a strictly increasing bijection from $(-m,m)$
onto $(0,+\infty)$.

On the other hand, the boundary condition associated with
$d_0^\gamma$ is
\begin{equation*}
    (\Gamma_1u,\Gamma_2u)
    \in
    \mathbb C
    \bigl(
        \cos\gamma,\,
        i\sin(\gamma+\pi/4)
    \bigr).
\end{equation*}
Consequently, if $\cos\gamma\neq0$, an eigenvalue in $(-m,m)$ exists
if and only if
\begin{equation*}
    \frac{\sin(\gamma+\pi/4)}{\cos\gamma}>0,
\end{equation*}
and in that case it is unique. For $\gamma\in[0,\pi)$ this is
equivalent to
\begin{equation*}
    \gamma\in[0,\pi/2)\cup(3\pi/4,\pi).
\end{equation*}
For
\begin{equation*}
    \gamma\in[\pi/2,3\pi/4]
\end{equation*}
there is no eigenvalue in the gap. At the endpoint
$\gamma=\pi/2$ one has $\Gamma_1u=0$, while at
$\gamma=3\pi/4$ one has $\Gamma_2u=0$; since the decaying solution
has both boundary coefficients non-zero for every $E\in(-m,m)$,
neither endpoint produces a bound state.

In particular, both distinguished extensions used in
Definition~\ref{defSA} have empty point spectrum in $(-m,m)$.
\end{remark}

%

\subsection{The relativistic Hankel transform}\label{sec:relhan}
We need to adapt the transform introduced in the predecessor of this paper to the massive setting. In order to simplify the presentation, we shall define the transform only for the distinguished self-adjoint extensions, as defined in Definition \ref{defSA}; also, without loss of generality we shall only discuss the case $\alpha\in(0,1/2]$, which corresponds to the choice $\gamma=3\pi/4$, the other one being completely analogous.

We recall that we are considering the following generalized eigenfunctions
for a given $\alpha\in(0,1)$ and $\rho\geq m$:
\begin{itemize}
  \item 
  for $k\ge1$,
  \begin{equation}\label{eq:psika}
    \Psi_{k}(\rho, r)=
    \frac{1}{\sqrt{2}}
    \begin{pmatrix}
     \sqrt{\rho-m} J_{k-\alpha}(\p r) \\
     - i\sqrt{\rho+m}J_{k+1-\alpha}(\p r)
    \end{pmatrix}
  \end{equation}
  \item 
  for $k\le-1$,
  \begin{equation}\label{eq:psikb}
    \Psi_{k}(\rho, r)=
    \frac{1}{\sqrt{2}}
    \begin{pmatrix}
     \sqrt{\rho-m} J_{\alpha-k}(\p r)\  \\
      i\sqrt{\rho+m}J_{\alpha-k-1}(\p r)
    \end{pmatrix}
  \end{equation}
  \item 
  for $k=0$,
  \begin{equation}\label{eq:psik0}
    \Psi_{0}(\rho, r)=
    \frac{1}{\sqrt{2}} 
    \begin{pmatrix}
  \sqrt{\rho-m}  J_{-\alpha}(\p r)\\
-  i \sqrt{\rho+m}    J_{1-\alpha}(\p r)
    \end{pmatrix}
  \end{equation}
\end{itemize}
where $\p=\p(\rho)=\sqrt{\rho^2-m^2}$\footnote{In what follows we shall systematically omit the dependence on $\rho$ in $\p$ to simplify the notation.}. For negative energies $\rho\leq -m$, the corresponding formulas are given by \eqref{neggen-eig'}-\eqref{negeq:psik0}.

Then we define the \emph{relativistic Hankel transform} as a
sequence of operators, $k\in \mathbb{Z}$
\begin{equation}\label{eq:Pkspa}
  \mathcal{P}_{k}:[L^{2}((0,\infty),rdr)]^{2}
    \to L^{2}(\sigma_{cont},d\rho)\footnote{In the massless convention of \cite[formula~(3.16)]{cacdanyinzha}, the
target measure is $|\rho|d\rho$. Here the factors $\sqrt{|\rho|\pm m}$
in the eigenfunctions absorb the spectral density, so the target measure is $d\rho$.}
\end{equation}
where $ \sigma_{cont}=]-\infty,-m]\cup[m,+\infty[$, 
acting on spinors
$f(r)=\big(\begin{smallmatrix}
   f_1 (r) \\ f_2 (r)\end{smallmatrix}\big)$, $r>0$,
as
\begin{equation}\label{hankel}
  \mathcal{P}_{k}f(\rho)=
  \begin{cases}
    \int_{0}^{+\infty}\Psi_{k}(\rho, r)^{*}\cdot f(r)rdr
     &\text{if $ \rho\geq m $,}\\
     \int_{0}^{+\infty}{\Psi_{k}}(\rho, r)^{*}
       \cdot f(r)rdr
     &\text{if $ \rho\leq -m $.}
  \end{cases}
\end{equation}
For $\rho>m$ and $k\ge1$, for example,
\begin{equation*}\begin{aligned}
\mathcal P_k f(\rho)
&=\frac1{\sqrt2}\int_0^\infty
\Bigl(\sqrt{\rho-m}J_{k-\alpha}(\p r)f_1(r)
+i\sqrt{\rho+m}J_{k+1-\alpha}(\p r)f_2(r)\Bigr)rdr\\
&=\frac1{\sqrt2}\Bigl(
\sqrt{\rho-m}(\mathcal H_{k-\alpha}f_1)(\p)
+i\sqrt{\rho+m}(\mathcal H_{k+1-\alpha}f_2)(\p)\Bigr),
\end{aligned}\end{equation*}
where $\mathcal{H}_{\nu}$ denotes the standard Hankel transform defined as
\begin{equation*}
  \mathcal H_\nu\phi(\xi)
  =\int_0^\infty J_\nu(r\xi)\phi(r)rdr,\qquad\xi>0.
\end{equation*}

Similar explicit formulas can be written for $\rho\leq -m$ and $k< 1$.

We now establish the following properties.
\begin{proposition}\label{pro:hankel}
  Let $\alpha\in(0,1/2]$ and $\gamma=\frac{3\pi}4$.
Let $\mathcal{P}_{k}$ be given by \eqref{hankel},
  with $\Psi_{k}$ given by \eqref{eq:psika}, \eqref{eq:psikb},
  \eqref{eq:psik0}, \eqref{neggen-eig'} and \eqref{negeq:psik0}

  Then the following properties hold:
  \begin{enumerate}
    \item 
    $\mathcal{P}_{k}:
      L^{2}((0,+\infty),rdr)^{2}\to
      L^{2}(\sigma_{cont},d\rho)$ 
    is a bounded bijection.

    \item 
    $\mathcal{P}_{k}$ is unitary.

    \item 
    We have
   $
    \mathcal P_k d_k^\gamma
    =
    \rho\,\mathcal P_k
    \qquad\text{on }D(d_k^\gamma).
$ \end{enumerate}
\end{proposition}

\begin{proof}
The presence of the mass term makes computations more delicate, so we prefer to include the full proof of this Proposition.
First of all, let us introduce the following convenient notation: we set 
\begin{equation}\label{hankpm}
(\mathcal H_\nu^\pm h)(\rho)
:=
\sqrt{\rho\pm m}\,
(\mathcal H_\nu h)\bigl(\sqrt{\rho^2-m^2}\bigr),
\qquad \rho>m.\end{equation}
Thus we can write for $k\ge0$, with $\nu=\nu(k)=k-\alpha$
\begin{equation}\label{eq:hankmat1}
  \mathcal{P}_{k} f(\rho)=
  2^{-\frac 12}
  \begin{cases}
    \mathcal{H}_{\nu}^-f_{1}(\rho)+i\mathcal{H}_{\nu+1}^+f_{2}(\rho)
      &\text{if $ \rho\geq m $,}\\
  \mathcal{H}_{\nu}^+f_{1}(|\rho|)-i\mathcal{H}_{\nu+1}^-f_{2}(|\rho|)
     &\text{if $ \rho\leq -m $}
  \end{cases}
\end{equation}
while for $k\le-1$, with $\nu=\alpha-k$, we have
\begin{equation}\label{eq:hankmat2}
  \mathcal{P}_{k} f(\rho)=
  2^{-\frac 12}
  \begin{cases}
    \mathcal{H}_{\nu}^-f_{1}(\rho)-i\mathcal{H}_{\nu-1}^+f_{2}(\rho)
      &\text{if $ \rho\geq m $,}\\
 \mathcal{H}_{\nu}^+f_{1}(|\rho|)+i\mathcal{H}_{\nu-1}^-f_{2}(|\rho|)
     &\text{if $ \rho\leq -m $}.
  \end{cases}
\end{equation}
Since each $\mathcal{H}_\nu$ is unitary on $L^2((0,\infty),rdr)$ for any $\nu\ge -\tfrac12$, we postpone the boundedness claim until after the isometry computation below. Indeed, that computation will show directly that
\begin{equation*}
\|\mathcal{P}_k f\|_{L^2(\sigma_{\mathrm{cont}},d\rho)}^2
=
\|f_1\|_{L^2((0,\infty),rdr)}^2+
\|f_2\|_{L^2((0,\infty),rdr)}^2,
\end{equation*}
and therefore $\mathcal{P}_k$ is bounded.

Let us now prove $\mathcal{P}_k$ is an isometry. Let us focus on the case $k\geq0$. From \eqref{eq:hankmat1}-\eqref{eq:hankmat2} we can write
\begin{equation}\label{explicitP}
\begin{pmatrix}
\mathcal{P}_kf(\rho)\\[1mm]
\mathcal{P}_kf(-\rho)
\end{pmatrix}
=
\frac1{\sqrt2}
M(\rho)
\begin{pmatrix}
(\mathcal{H}_{\nu}f_1)(\sqrt{\rho^2-m^2})\\[1mm]
(\mathcal{H}_{\nu+1}f_2)(\sqrt{\rho^2-m^2})
\end{pmatrix}
\end{equation}
with $\nu=k-\alpha$ and where we have defined
\begin{equation*}
M(\rho)=
\begin{pmatrix}
\sqrt{\rho-m} & i\sqrt{\rho+m}\\
\sqrt{\rho+m} & -i\sqrt{\rho-m}
\end{pmatrix}.
\end{equation*}
A direct computation gives
\begin{equation}\label{Mmatrix}
\frac12M(\rho)^*M(\rho)=\rho\,I_2.
\end{equation}
Hence
\begin{equation}\label{isometry}
|\mathcal{P}_kf(\rho)|^2+|\mathcal{P}_kf(-\rho)|^2
=
\rho\Big(|(\mathcal H_\nu f_1)(\p)|^2+|(\mathcal H_{\nu+1}f_2)(\p)|^2\Big).
\end{equation}
The change of variables
\begin{equation*}
\p=\sqrt{\rho^2-m^2},\qquad \rho=\sqrt{\p^2+m^2}
\end{equation*}
gives
\begin{equation*}
d\rho=\frac{\p}{\rho}\,d\p
\end{equation*}
and thus from \eqref{isometry} we get
\begin{equation*}\begin{aligned}
\int_m^\infty \Big(|\mathcal{P}_kf(\rho)|^2+|\mathcal{P}_kf(-\rho)|^2\Big)\,d\rho
&=
\int_m^\infty \rho\Big(|(\mathcal H_\nu f_1)(\p)|^2+|(\mathcal H_{\nu+1}f_2)(\p)|^2\Big)\,d\rho
\\
&=
\int_0^\infty \Big(|\mathcal{H}_\nu f_1(\p)|^2+|\mathcal{H}_{\nu+1}f_2(\p)|^2\Big)\,\p\,d\p.
\end{aligned}\end{equation*}
Since the Hankel transforms $\mathcal{H}_\mu$ are unitary on $L^2((0,\infty),r\,dr)$ for
every $\mu\ge -\tfrac12$, and here $\nu,\nu+1\ge -\tfrac12$ because
$\alpha\in(0,\tfrac12]$, it follows that
\begin{equation*}
\|\mathcal{P}_kf\|_{L^2(\sigma_{\mathrm{cont}},d\rho)}^2
=
\|f_1\|_{L^2((0,\infty),r\,dr)}^2
+
\|f_2\|_{L^2((0,\infty),r\,dr)}^2.
\end{equation*}
Thus $\mathcal{P}_k$ is an isometry for every $k\ge 0$.

The same computation applies for $k\le-1$, with the signs in
\eqref{eq:hankmat2}. In fact, $\mathcal P_k$ is onto: the two Hankel transforms
are unitary, the change of variables $\p\,d\p=\rho\,d\rho$ is unitary after
the factor $\sqrt\rho$ is included, and $M(\rho)/\sqrt{2\rho}$ is a unitary
$2\times2$ matrix. Thus their composition is unitary onto the two spectral branches.

For clarity, put $\phi_\pm(\rho)=\phi(\pm\rho)$ for $\rho>m$. For $k\ge0$,
with $\nu=k-\alpha$, the inverse is

\begin{equation}\label{eq:inverse-hankel-positive}
\begin{aligned}
(\mathcal P_k^{-1}\phi)_1(r)
&=\frac1{\sqrt2}\int_m^\infty J_\nu(\p r)
\Bigl(\sqrt{\rho-m}\phi_+(\rho)+\sqrt{\rho+m}\phi_-(\rho)\Bigr)d\rho,\\
(\mathcal P_k^{-1}\phi)_2(r)
&=\frac1{\sqrt2}\int_m^\infty J_{\nu+1}(\p r)
\Bigl(-i\sqrt{\rho+m}\phi_+(\rho)+i\sqrt{\rho-m}\phi_-(\rho)\Bigr)d\rho.
\end{aligned}
\end{equation}

For $k\le-1$, set $\nu=\alpha-k$; the first line is unchanged, and the
second becomes

\begin{equation*}
(\mathcal P_k^{-1}\phi)_2(r)
=\frac1{\sqrt2}\int_m^\infty J_{\nu-1}(\p r)
\Bigl(i\sqrt{\rho+m}\phi_+(\rho)-i\sqrt{\rho-m}\phi_-(\rho)\Bigr)d\rho.
\end{equation*}

These formulas are first understood for compactly supported spectral data
and then by $L^2$ continuity.

It remains to prove the intertwining identity. For $k\ne0$, work first on
the core $[C_c^\infty((0,\infty))]^2$. For $k=0$, use this core together
with the additional domain vector $(K_\alpha(\sqrt2 mr),0)^T$.
Green's identity gives the boundary term
\begin{equation*}
i\Bigl[r\bigl(\overline{\Psi_{k,1}(\rho,r)}f_2(r)
+\overline{\Psi_{k,2}(\rho,r)}f_1(r)\bigr)\Bigr]_{0}^{\infty}.
\end{equation*}
It vanishes for the compactly supported core functions. For the additional
critical vector, its only nonzero contribution is a constant multiple of
$rJ_{1-\alpha}(\p r)K_\alpha(\sqrt2 mr)$; this is
$O(r^{2-2\alpha})$ at zero and decays at infinity. Consequently,

\begin{equation*}
(\mathcal P_kd_k^\gamma f)(\rho)
=\int_0^\infty(d_k^\gamma\Psi_k(\rho,r))^*f(r)rdr
=\rho\,(\mathcal P_k f)(\rho)
\end{equation*}
on the core. The equality extends in the graph norm to 
$D(d_k^\gamma)$,
because $\mathcal P_k$ is unitary and multiplication by $\rho$ is closed.
Finally, both operators are self-adjoint, so this inclusion is an equality
of operators and their domains. This proves the proposition.
\end{proof}

\subsection{Representation of the propagator.}

Let us assume for the rest of this subsection that
$\alpha\in(0,1/2]$ and let us choose
$\gamma=\frac{3\pi}{4}$ in Theorem~\ref{the:sadjextDAB}. We shall simply denote the operator $\mathcal{D}_A$ in what follows. For the sign convention used in the estimates, let $u$ solve
\begin{equation}\label{eq:Dirac}
\begin{cases}
i\partial_t u+\mathcal{D}_{A} u=0,\\
u|_{t=0}=f(x)=(\phi, \psi)^T \in D(\mathcal{D}_{A})\subset (L^2(\R^2))^2
\end{cases}
\end{equation}
where the domain is described in Theorem~\ref{the:sadjextDAB}.

We shall need to treat separately
the case $k\neq0$, which we call
the {\em noncritical part}, and can be handled by relying on the existing theory for the Klein--Gordon flow, from the {\em critical part}, that corresponds to $k=0$, in which the singularity of the domain comes into play: this sector requires a more careful, ``ad hoc" analysis.

To separate the noncritical part we apply
the orthogonal projector $P_{\bot}$  defined in \eqref{def-pro},
which cuts away the spherical harmonic spinors of order $0$
\begin{equation*}
  f_{\bot}=P_{\bot}f=  \sum_{k\in \mathbb{Z}\setminus\{0\}}
  \begin{pmatrix}
    e^{ik \theta}\ \phi_{k}(r) \\
    e^{i(k+1)\theta}\ \psi_{k}(r).
  \end{pmatrix}
\end{equation*} 
This corresponds to the choice of data
\begin{equation*}
  f_0=P_{0}f=
  \begin{pmatrix}
  \phi_{0}(r) \\
    e^{i\theta}\ \psi_{0}(r)
  \end{pmatrix}=
  c\sqrt m
    \begin{pmatrix}
      \cos \gamma \cdot K_{\alpha}((\sqrt2 mr) ) \\
    0
    \end{pmatrix}
  +\begin{pmatrix}
  \tilde\phi_{0}(r) \\
    e^{i\theta}\ \tilde\psi_{0}(r)
  \end{pmatrix}.
\end{equation*}
By orthogonality, the full flow is the sum
of the noncritical and the critical components:
 \begin{equation*}
e^{it\mathcal{D}_{A}} f=e^{it\mathcal{D}_{A}}P_0 f+e^{it\mathcal{D}_{A}}P_{\bot} f=e^{it\mathcal{D}_{A}} f_0+e^{it\mathcal{D}_{A}} f_{\bot}.
  \end{equation*}

We start by providing an explicit representation of the propagator $e^{it\mathcal{D}_{A}}$ via the generalized eigenfunctions.
In the following statement we denote polar coordinates on
$\mathbb{R}^{2}$ by $(r,\theta)$, like in
$x=r(\cos\theta,\sin\theta)$ and 
$y=r_2(\cos\theta_2,\sin\theta_2)$.

\begin{proposition}\label{prop:Dp} 
Let $\alpha\in(0, 1/2]$ and $\gamma=\frac{3\pi}4$.
Let $\Phi_{k}(\theta)$ be the harmonics
\begin{equation}\label{eq:defPhik}
  \Phi_k(\theta)=\frac{1}{\sqrt{2\pi}}
  \left(\begin{matrix} e^{ik\theta}&0\\ 
  0&e^{i(k+1)\theta}\end{matrix}\right).
\end{equation}
Then the solution $u(t,x)$ of
\eqref{eq:Dirac} can be represented as 
\begin{equation}\label{D-propa}
u(t, x)=e^{it\mathcal{D}_{A}}f=\int_{0}^\infty \int_{0}^{2\pi}{\bf K}(t,r,\theta,r_2,\theta_2) f(r_2,\theta_2)  d\theta_2 \;r_2dr_2,
\end{equation}
\begin{equation}\label{D-kernel}
{\bf K}(t,r,\theta,r_2,\theta_2) =\sum_{k\in\Z} 
\Phi_k(\theta) {\bf K}_{k}(t,r,r_2)\overline{\Phi_k(\theta_2)},
\end{equation}
where, setting
\begin{equation}\label{eq:defdiagmassive}
  d_{\nu}^{\pm}(t,r,r_2)
  :=
  \frac12\int_m^{\infty}
  \Big((\rho\mp m)e^{it\rho}+(\rho\pm m)e^{-it\rho}\Big)
  J_{\nu}(\p r)J_{\nu}(\p r_2)\,d\rho,
\end{equation}
\begin{equation}\label{eq:defoffdiagmassive}
  m_{\nu,\mu}^{s}(t,r,r_2)
  :=
  -\frac{i}{2}\int_m^{\infty}
  \p\big(e^{it\rho}-e^{-it\rho}\big)
  J_{\nu}(\p r)J_{\mu}(\p r_2)\,d\rho,
\end{equation}
with $\p=\sqrt{\rho^2-m^2}$, one has the following explicit formulas.

\smallskip
\noindent
If $k\ge 1$, setting $\nu=k-\alpha$, then
\begin{equation}\label{kernel'}
{\bf K}_{k}(t, r,r_2)=
\begin{pmatrix}
  d_{\nu}^{+}(t,r,r_2) & -\,m_{\nu,\nu+1}^{s}(t,r,r_2)  \\
  m_{\nu+1,\nu}^{s}(t,r,r_2) &  d_{\nu+1}^{-}(t,r,r_2)
\end{pmatrix}.
\end{equation}
If $k\le -1$, setting $\nu=\alpha-k$, then
\begin{equation}\label{eq:kernel-neg}
{\bf K}_{k}(t, r,r_2)=
\begin{pmatrix}
  d_{\nu}^{+}(t,r,r_2) & m_{\nu,\nu-1}^{s}(t,r,r_2)  \\
  -\,m_{\nu-1,\nu}^{s}(t,r,r_2) &  d_{\nu-1}^{-}(t,r,r_2)
\end{pmatrix}.
\end{equation}
For $k=0$, one has
\begin{equation}\label{kernel'0}
  {\bf K}_{0}(t, r,r_2)=
  \begin{pmatrix}
    d_{-\alpha}^{+}(t,r,r_2) &-\,m_{-\alpha,1-\alpha}^{s}(t,r,r_2)\\
    m_{1-\alpha,-\alpha}^{s}(t,r,r_2) &  d_{1-\alpha}^{-}(t,r,r_2)
  \end{pmatrix}.
\end{equation}
\end{proposition}

\begin{proof}
By Proposition~\ref{pro:hankel}, the radial kernel is
\begin{equation*}
\mathbf K_k(t,r,r_2)
=\int_m^\infty\left(
e^{it\rho}\Psi_k(\rho,r)\Psi_k(\rho,r_2)^*
+e^{-it\rho}\Psi_k(-\rho,r)\Psi_k(-\rho,r_2)^*
\right)d\rho.
\end{equation*}

For $k\ge0$, the upper-right entry is
$\frac{i}{2}\int_m^\infty\p(e^{it\rho}-e^{-it\rho})
J_\nu(\p r)J_{\nu+1}(\p r_2)d\rho=-m_{\nu,\nu+1}^s$;
the lower-left entry has the opposite sign and the Bessel orders reversed.
The two diagonal entries are $d_\nu^+$ and $d_{\nu+1}^-$.
For $k\le-1$, the signs of the off diagonal entries are reversed, as in
\eqref{eq:kernel-neg}. Reconstructing the angular modes yields
\eqref{D-kernel}. At $t=0$, the off diagonal terms vanish and the diagonal
ones reduce to Hankel inversion, which also checks the normalization.
\end{proof}

\begin{remark}
Notice the equivalent representations
\begin{equation*}
  d_{\nu}^{\pm}(t,r,r_2)
  =
  \int_m^{\infty}
  \big(\rho\cos(t\rho)\mp i m\sin(t\rho)\big)
  J_{\nu}(\p r)J_{\nu}(\p r_2)\,d\rho,
\end{equation*}
\begin{equation*}
  m_{\nu,\mu}^{s}(t,r,r_2)
  =
  \int_m^{\infty}
  \p\sin(t\rho)
  J_{\nu}(\p r)J_{\mu}(\p r_2)\,d\rho.
\end{equation*}
\end{remark}

For the spacetime estimates it is useful to record an equivalent
single-Bessel representation. Let $\mathbf f_0(r)$ be the radial coefficient
of $P_0f$ in the normalized angular basis $\Phi_0$, and set

\begin{equation*}
\omega(\xi)=\sqrt{m^2+\xi^2},\qquad
g_\pm(\xi)=\omega(\xi)^{-1/2}
(\mathcal P_0\mathbf f_0)(\pm\omega(\xi)).
\end{equation*}
Then
\begin{equation}\label{eq:momentum-unitarity}
\|g_+\|_{L^2(\xi d\xi)}^2+\|g_-\|_{L^2(\xi d\xi)}^2
=\|P_0f\|_{[L^2(\R^2)]^2}^2.
\end{equation}
The radial coefficient of the evolved spinor is
\begin{equation}\label{eq:single-bessel-flow}
\mathbf u_0(t,r)=\sum_{\sigma\in\{+,-\}}
\int_0^\infty e^{i\sigma t\omega(\xi)}
\binom{b_{\sigma,1}(\xi)J_{-\alpha}(r\xi)}
{b_{\sigma,2}(\xi)J_{1-\alpha}(r\xi)}g_\sigma(\xi)\xi\,d\xi,
\end{equation}
where
\begin{equation*}
\begin{aligned}
b_{+,1}&=\sqrt{\frac{\omega-m}{2\omega}},&
b_{+,2}&=-i\sqrt{\frac{\omega+m}{2\omega}},\\
b_{-,1}&=\sqrt{\frac{\omega+m}{2\omega}},&
b_{-,2}&=i\sqrt{\frac{\omega-m}{2\omega}}.
\end{aligned}
\end{equation*}
These coefficients are bounded, with uniformly bounded derivatives after
rescaling on either high frequency momentum shells or threshold momentum
shells. Thus the off diagonal terms introduce no additional $L^2$ input
norm in the spacetime estimates.

\begin{remark}\label{Generapspecrep}
For a generic self-adjoint extension of the Dirac operator, the propagator has the same form for the noncritical component $k\neq0$. In the singular channel $k=0$ on the other hand, the continuous kernel must be built from the generalized eigenfunctions satisfying the boundary condition \eqref{eq:critical-bc}, and therefore depends on the chosen extension parameter. Moreover, when a bound state is present, one has to add the corresponding point-spectrum contribution, namely
\begin{equation*}
(e^{it\mathcal D_A}f)(x)
=\int_0^\infty\int_0^{2\pi}{\bf K}(t,r,\theta,r_2,\theta_2)f(r_2,\theta_2)\,d\theta_2\,r_2dr_2
+e^{it\lambda}\langle\phi_\lambda,f\rangle\phi_\lambda(x),
\end{equation*}
where $\lambda$ is the unique eigenvalue (for those extensions for which it exists) and $\phi_\lambda$ is the associated eigenfunction. Clearly the additional term, whenever present, will produce a non-dispersive component of the solution.
\end{remark}

\section{Localized dispersive estimates}\label{sec:proofdisp}

In this section we adapt Section~5 of the massless paper to the massive setting.
Throughout we work with the distinguished self-adjoint extension fixed in
Definition~\ref{defSA}; in particular we assume
\begin{equation*}
\alpha\in(0,\tfrac12],
\qquad
\gamma=\frac{3\pi}{4}
\end{equation*}
and we shall simply denote $\mathcal{D}_A=\mathcal{D}_{A,\gamma}$ for simplicity.
 The case $1/2<\alpha<1$ can be obtained by
replacing $\alpha$ with $1-\alpha$, interchanging the spinor components,
and interchanging the positive and negative energy branches (equivalently,
reversing time in the radial representation). More explicitly, if
$d_{0,\alpha}$ denotes the critical radial operator at flux $\alpha$ and
$S=\bigl(\begin{smallmatrix}0&1\\-1&0\end{smallmatrix}\bigr)$, then
\begin{equation*}
d_{0,\alpha}=-S d_{0,1-\alpha}S^{-1}.
\end{equation*}
The matrix $S$ maps the distinguished radial domains onto one another,
and the angular reconstruction only adds factors of modulus one. Thus
the reflection preserves the norms and interchanges the two diagonal
weights. As in the massless case, we separate the regular sector $P_{\bot}$ and the
singular sector $P_{0}$. The main two differences are the following:

\begin{itemize}
\item On the noncritical (regular) sector, the Dirac flow must be split into its positive and
negative spectral branches, exactly as in the erratum to the massless paper. Then, dispersive estimates can be obtained by applying the results from \cite{GaoYinZhangZheng} for the Klein--Gordon flow.
\item On the radial (singular) sector, some explicit bounds on the generalized eigenfunctions need to be derived. \end{itemize}

\subsection{The regular sector: proof of \eqref{eq:new121}}\label{subsec:regular-dispersive}

We denote by $\Pi_{\bot}^{\pm}$ the spectral projections of $\mathcal D_A$ on
$P_{\bot}[L^2(\R^2)]^2$ corresponding to the positive and negative spectrum,
respectively, and set
\begin{equation*}
\mathcal H_\perp:=\Ran P_\perp,
\qquad
D_\perp:=\mathcal D_A\big|_{\mathcal H_\perp}.
\end{equation*}

Let $H_A^{\rm F}$ denote the Friedrichs realization of the scalar
Aharonov--Bohm Hamiltonian and define, on the full spinor space,
\begin{equation*}
\mathscr L_m
:=
\begin{pmatrix}
H_A^{\rm F}+m^2&0\\
0&H_A^{\rm F}+m^2
\end{pmatrix}.
\end{equation*}
On the regular subspace we set
\begin{equation*}
L_m:=\mathscr L_m\big|_{\mathcal H_\perp}.
\end{equation*}
By the partial-wave decomposition,
\begin{equation*}
D_\perp^2=L_m,
\qquad
|D_\perp|=\sqrt{L_m}.
\end{equation*}
Equivalently,
\begin{equation*}
L_m
=
\bigoplus_{k\neq 0}
\begin{pmatrix}
H_{|k-\alpha|}+m^2 & 0\\
0 & H_{|k+1-\alpha|}+m^2
\end{pmatrix}.
\end{equation*}
Hence, by the spectral theorem,
\begin{equation}\label{eq:regular-splitting-massive}
e^{it\mathcal D_A}P_{\bot}
=
e^{it\sqrt{\mathscr L_m}}\Pi_{\bot}^{+}
+
e^{-it\sqrt{\mathscr L_m}}\Pi_{\bot}^{-}.
\end{equation}
Here the half-wave groups generated by $\sqrt{\mathscr L_m}$ may be used because
their restriction to $\mathcal H_\perp$ is generated by
$|D_\perp|=\sqrt{L_m}$. This is the exact analogue of the corrected
regular-sector splitting in the massless erratum.

Then, the time decay estimate on the regular sector, i.e. estimate \eqref{eq:new121} is a consequence of the following
\begin{proposition}\label{prop:regular-dispersive-massive}
Let $\Pi_\perp^\pm$ be the positive and negative regular spectral projections. 
Fix $\phi_0,\widetilde\phi_0\in C_c^\infty([0,\infty))$ such that
\begin{equation*}
0\le \phi_0,\widetilde\phi_0\le 1,\qquad
\widetilde\phi_0\,\phi_0=\phi_0,
\end{equation*}
\begin{equation*}
\phi_0(\lambda)=1\quad\text{for }0\le \lambda\le 4m,
\qquad
\supp \phi_0\subset [0,8m],
\end{equation*}
and let $\{\varphi_j\}_{j\in\mathbb Z}$, $\{\widetilde\varphi_j\}_{j\in\mathbb Z}$ be the
dyadic cutoffs as in Theorem \ref{thm:disper}.

Then there exists a constant $C>0$ such that, for every $t\in\R$, the following estimates hold.

\smallskip
\noindent
{\rm(i) Low frequencies.} For every
$f\in [L^1(\R^2)]^2\cap [L^2(\R^2)]^2$,
\begin{equation}\label{eq:regular-low-dispersive-massive}
\big\|
\phi_0(|\mathcal D_A|)\,e^{it\mathcal D_A}\,\Pi_\perp^\pm f
\big\|_{[L^\infty(\R^2)]^2}
\le
C(1+|t|)^{-1}
\big\|
\widetilde\phi_0(|\mathcal D_A|)\,\Pi_\perp^\pm f
\big\|_{[L^1(\R^2)]^2}.
\end{equation}

\smallskip
\noindent
{\rm(ii) High frequencies.} For every $j\in\Z$ such that $2^j\ge 4m$ and every
$f\in [L^1(\R^2)]^2\cap [L^2(\R^2)]^2$,
\begin{equation}\label{eq:regular-high-dispersive-massive}
\big\|
\varphi_j(|\mathcal D_A|)\,e^{it\mathcal D_A}\,\Pi_\perp^\pm f
\big\|_{[L^\infty(\R^2)]^2}
\le
C\,2^{2j}(1+2^j|t|)^{-1/2}
\big\|
\widetilde\varphi_j(|\mathcal D_A|)\,\Pi_\perp^\pm f
\big\|_{[L^1(\R^2)]^2}.
\end{equation}
\end{proposition}

\begin{proof}
We apply the scalar kernel estimates underlying
Propositions~3.3 and~3.4 of \cite{GaoYinZhangZheng}. Those propositions use
cutoffs in the momentum $\sqrt{H_A^{\rm F}}$, whereas the present cutoffs
are in the energy $\sqrt{H_A^{\rm F}+m^2}$. At low energy the composed
symbol $\phi_0(\sqrt{\xi^2+m^2})$ is a smooth even function of $\xi$.
For $2^j\ge4m$, its dyadic analogue is supported where $\xi\simeq2^j$
and has uniformly bounded rescaled derivatives. Thus the same kernel proofs
apply; one is not identifying the two cutoff variables. Indeed, by the regular/singular decomposition established above, on the regular
subspace the flow is reduced to the purely magnetic Klein--Gordon propagator on the
noncritical angular channels, and the two regular branches are given by
\begin{equation*}
e^{it\mathcal D_A}\Pi_\perp^\pm
=
e^{\pm it\sqrt{L_m}}\Pi_\perp^\pm .
\end{equation*}
Thus, after angular decomposition, the regular part is treated componentwise by the
same scalar kernels as in the Klein--Gordon case, and we can thus conclude that:
\begin{itemize}
\item 
For the low-frequency part, \eqref{eq:regular-low-dispersive-massive} follows from the
low energy dispersive estimate of Proposition~3.3 in \cite{GaoYinZhangZheng}, since the
single cutoff $\phi_0$ absorbs all finitely many shells intersecting the threshold region
$|\lambda|\sim m$.
\item
For the high frequency part, if $2^j\ge 4m$, then
\eqref{eq:regular-high-dispersive-massive} follows from Proposition~3.4 in \cite{GaoYinZhangZheng} with $\eta=0$, which gives, as decay factor, 
\begin{equation*}
2^{\frac{3j}{2}}\bigl(2^{-j}+|t|\bigr)^{-1/2}
=
2^{2j}(1+2^j|t|)^{-1/2}.
\end{equation*}
\end{itemize}
Applying these scalar estimates to each component of the regular flow and summing the
two components gives \eqref{eq:regular-low-dispersive-massive} and
\eqref{eq:regular-high-dispersive-massive}.
The finitely many nonzero dyadic shells with $2^j<4m$ are treated by the
same low energy kernel estimate with their own composed cutoffs.
Since $2^j\simeq_m1$ there, their $(1+|t|)^{-1}$ bound implies
\eqref{eq:new121}; shells with $2^{j+1}\le m$ vanish.
\end{proof}
%

\subsection{Weighted estimate in the singular sector: proof of \eqref{est:dis0}}\label{subsec:weighted-dispersive}

\medskip
\noindent
Before entering the proof, let us briefly comment on the differences with respect to the
massless case. 
The genuinely new feature of the
massive case is that the oscillatory phase is no longer $e^{\pm it\rho}$ with Bessel
argument $\rho r$, but rather
\begin{equation*}
e^{\pm it\sqrt{\p^2+m^2}},
\qquad \p=\p(\rho)=\sqrt{\rho^2-m^2},
\end{equation*}
and the Bessel functions depend on $\p r$. Thus the spectral variable of the oscillation and
the natural Bessel variable no longer coincide; in other words, in the massive case we are losing the natural scaling structure. As it is often the case for the analysis of Klein--Gordon or massive Dirac equations, the dynamics for high and low frequencies behave quite differently.  For this reason, unlike in the
massless case, we find it convenient to split the proof into two parts:
\begin{itemize}
\item a finite family of ``threshold energies"
$2^j\simeq m$, treated separately,
\item the genuine high frequency region $2^j\gg m$,
where the argument reduces to a smooth perturbation of the massless one after passing to
the momentum variable $\p$.
\end{itemize}
\medskip

For the singular sector we use directly the kernel representation of
Proposition~\ref{prop:Dp}. Throughout this subsection we assume again that 
\begin{equation*}
\alpha\in(0,\tfrac12],
\qquad
\gamma=\frac{3\pi}{4}
\end{equation*}
and we shall simply denote $\mathcal{D}_A=\mathcal{D}_{A,\gamma}$ for simplicity.
The case $\alpha\in(1/2,1)$ can be dealt with symmetrically.
We set
\begin{equation}\label{eq:Wj-massive}
W_{j}(r):=
\begin{pmatrix}
(1+(2^{j}r)^{-\alpha})^{-1} & 0\\
0 & 1
\end{pmatrix},
\qquad
W_{j}(r)^{-1}=
\begin{pmatrix}
1+(2^{j}r)^{-\alpha} & 0\\
0 & 1
\end{pmatrix}.
\end{equation}
Let
\begin{equation*}
f_{j}=\widetilde\varphi_{j}(|\mathcal{D}_{A}|)P_{0}f=
\begin{pmatrix}
\phi_{j}(r)\\
e^{i\theta}\psi_{j}(r)
\end{pmatrix}.
\end{equation*}
Using \eqref{kernel'0} we can write
\begin{equation*}
\varphi_{j}(|\mathcal{D}_{A}|)e^{it\mathcal{D}_{A}}P_{0}f
=
\int_{0}^{\infty}\int_{0}^{2\pi}
\Phi_{0}(\theta){\bf K}_{0,m}^{(j)}(t,r,r_{2})\overline{\Phi_{0}(\theta_{2})}
 f_{j}(r_{2},\theta_{2})\,d\theta_{2}\,r_{2}dr_{2},
\end{equation*}
where
\begin{equation}\label{eq:singular-localized-kernel-massive}
\Phi_{0}(\theta){\bf K}_{0,m}^{(j)}(t,r,r_{2})\overline{\Phi_{0}(\theta_{2})}
=
\begin{pmatrix}
F_{\alpha,j}^{m}(t;r,r_{2}) & -e^{-i\theta_{2}}S_{\alpha,j}^{m}(t;r,r_{2})\\
e^{i\theta}\widetilde S_{\alpha,j}^{m}(t;r,r_{2}) & e^{i(\theta-\theta_{2})}E_{\alpha,j}^{m}(t;r,r_{2})
\end{pmatrix}
\end{equation}
with
\begin{equation}\label{eq:Fj-massive}
F_{\alpha,j}^{m}(t;r,r_{2})=
\frac1{4\pi}
\int_{m}^{\infty}
\Big((\rho-m)e^{it\rho}+(\rho+m)e^{-it\rho}\Big)
J_{-\alpha}(\p r)J_{-\alpha}(\p r_{2})\varphi_{j}(\rho)\,d\rho,
\end{equation}
\begin{equation}\label{eq:Ej-massive}
E_{\alpha,j}^{m}(t;r,r_{2})=
\frac1{4\pi}
\int_{m}^{\infty}
\Big((\rho+m)e^{it\rho}+(\rho-m)e^{-it\rho}\Big) 
J_{1-\alpha}(\p r)J_{1-\alpha}(\p r_{2})\varphi_{j}(\rho)\,d\rho,
\end{equation}
\begin{equation}\label{eq:Sj-massive}
S_{\alpha,j}^{m}(t;r,r_{2})=
\frac1{2\pi}
\int_{m}^{\infty}
\p\sin(t\rho)
J_{-\alpha}(\p r)J_{1-\alpha}(\p r_{2})\varphi_{j}(\rho)\,d\rho,
\end{equation}
\begin{equation}\label{eq:Stildej-massive}
\widetilde S_{\alpha,j}^{m}(t;r,r_{2})=
\frac1{2\pi}
\int_{m}^{\infty}
\p\sin(t\rho)
J_{1-\alpha}(\p r)J_{-\alpha}(\p r_{2})\varphi_{j}(\rho)\,d\rho.
\end{equation}

We shall use the standard identities and asymptotic expansions for
Bessel functions; see \cite{Watson}. The following bounds are
recorded in the form used in \cite[Lemma~5.1]{cacdanyinzha}.
\begin{lemma}
\label{lem:Bessel-massive-sing}
Let $\nu_-:=\max\{0,-\nu\}$. For every fixed $\nu\in\mathbb R$ and every
$x>0$, one has
\begin{equation}\label{eq:Bessel-basic-massive}
        |J_\nu(x)|
        \le C_\nu x^\nu\langle x\rangle^{-\nu-\frac12},
        \qquad
        |J_\nu(x)|
        \le C_\nu(1+x^{-\nu_-}),
\end{equation}
and
\begin{equation}\label{eq:Bessel-derivative-massive}
        |J_\nu'(x)|
        \le
        C_\nu x^{\nu-1}\langle x\rangle^{-\nu+\frac12}.
\end{equation}
Moreover, for $x\ge1$, one can write
\begin{equation}\label{eq:Bessel-osc-massive}
        J_\nu(x)
        =
        x^{-1/2}\bigl(e^{ix}a_+(x)+e^{-ix}a_-(x)\bigr),
\end{equation}
where $a_\pm$ are smooth and satisfy, for every integer $N\ge1$,
\begin{equation*}
        |a_\pm(x)|\le C_\nu,
        \qquad
        |\partial_x^N a_\pm(x)|\le C_{\nu,N}x^{-N-1},
        \qquad x\ge1 .
\end{equation*}
\end{lemma}

\begin{proof}
This is exactly Lemma~5.1 of \cite{cacdanyinzha}; we only record it here for
convenience.
\end{proof}

Next, we need the following

\begin{lemma}[Fixed-shell Bessel estimate]
\label{lem:fixed-shell-bessel-massive}
Let $\nu,\mu\ge -\frac12$, and set
\begin{equation*}
        \nu_-:=\max\{0,-\nu\},\qquad
        \mu_-:=\max\{0,-\mu\}.
\end{equation*}
Let $b$ be a smooth function supported in a fixed compact interval
$I\Subset(0,\infty)$. Let $\omega$ be real and smooth on a neighbourhood
of $I$, with
\begin{equation*}
        |\omega'(\lambda)|\ge c>0,
        \qquad
        |\partial_\lambda^N\omega(\lambda)|\le C_N,
        \qquad N\ge2 .
\end{equation*}
Then, uniformly for $T\ge0$ and $R_1,R_2>0$,
\begin{equation*}
\left|
\int_0^\infty
e^{iT\omega(\lambda)}
b(\lambda)
J_\nu(R_1\lambda)J_\mu(R_2\lambda)\lambda\,d\lambda
\right|
\le
C(1+T)^{-1/2}
(1+R_1^{-\nu_-})(1+R_2^{-\mu_-}).
\end{equation*}
The constant $C$ depends only on $\nu,\mu$, on the fixed compact support of $b$,
and on finitely many seminorms of $b$ and $\omega$. In particular, the estimate is
uniform for families $b_\kappa,\omega_\kappa$ satisfying the above bounds uniformly.
\end{lemma}

\begin{proof}
This is the fixed-frequency estimate of Lemma~5.2 from \cite{cacdanyinzha}, with two
minor harmless modifications: the two Bessel orders may be different, and the linear
phase is replaced by a smooth fixed-shell phase$\omega(\lambda)$ with nonvanishing first derivative.

For completeness, we include a sketch of the argument. 

If $T\le1$, the estimate follows directly
from
\begin{equation*}
        |J_\sigma(z)|\lesssim 1+z^{-\sigma_-}.
\end{equation*}
Assume $T\ge1$. We distinguish two cases. If $R_1+R_2\gtrsim T$, then at least
one Bessel factor has large argument of size $\gtrsim T$, unless $T\lesssim1$, which
has already been treated. Using the large-argument expansion \eqref{eq:Bessel-osc-massive} 
for that factor, and the pointwise bound for the other one, gives the desired
$T^{-1/2}$ decay.

If instead $R_1+R_2\ll T$, we use the small-argument bounds for the Bessel factors
with $R_i\lesssim1$, and the large-argument expansion for those with $R_i\gtrsim1$.
Each resulting phase has the form
\begin{equation*}
        \Phi(\lambda)
        =
        T\omega(\lambda)+\varepsilon_1R_1\lambda+\varepsilon_2R_2\lambda,
        \qquad
        \varepsilon_i\in\{-1,0,1\}.
\end{equation*}
Choosing the implicit constant in $R_1+R_2\ll T$ sufficiently small, we have
\begin{equation*}
        |\Phi'(\lambda)|
        \ge
        T|\omega'(\lambda)|-R_1-R_2
        \gtrsim T
\end{equation*}
on $\operatorname{supp}b$. One integration by parts then gives the stronger bound
\begin{equation*}
        O(T^{-1})
        (1+R_1^{-\nu_-})(1+R_2^{-\mu_-}).
\end{equation*}
This proves the estimate.
\end{proof}
We are now in a position to prove the key estimates for the kernel functions; we isolate the contribution of the finitely many dyadic shells meeting the threshold
$\rho=m$ first.

\begin{lemma}[Threshold-energies bound]
\label{lem:threshold-massive-sing}
Let $\nu,\mu\geq-\tfrac12$, let $C_m>0$, and let $\chi\in C_c^\infty([0,\infty))$ with $\supp\chi\subset[0,C_m]$. Define
\begin{equation*}
I_{\nu,\mu}^{\pm}(t;r_1,r_2)
:=
\int_0^\infty
e^{\pm it\sqrt{\p^2+m^2}}
J_\nu(\p r_1)J_\mu(\p r_2)\chi(\p)\,\p\,d\p .
\end{equation*}
Then there exists a constant $C=C(m,\nu,\mu,\chi)$ such that
\begin{equation*}
|I_{\nu,\mu}^{\pm}(t;r_1,r_2)|
\le
C(1+|t|)^{-1/2}
\bigl(1+r_1^{-\nu_-}\bigr)\bigl(1+r_2^{-\mu_-}\bigr),
\end{equation*}
where $\nu_-:=\max\{0,-\nu\},\ \mu_-:=\max\{0,-\mu\}$.
\end{lemma}

\begin{proof}
The case $|t|\le1$ follows directly from Lemma~\ref{lem:Bessel-massive-sing}. Indeed,
\begin{equation*}
\begin{aligned}
|I^\pm_{\nu,\mu}(t;r_1,r_2)|
&\le
C\int_0^{C_m}
\p(1+(\p r_1)^{-\nu_-})(1+(\p r_2)^{-\mu_-})\,d\p  \\
&\le
C(1+r_1^{-\nu_-})(1+r_2^{-\mu_-}),
\end{aligned}
\end{equation*}
because $\nu_-+\mu_-\le1<2$, so the integrand is locally integrable
at $\p=0$.

We now assume $|t|\ge1$. First suppose that
\begin{equation*}
        \nu_-+\mu_-<1.
\end{equation*}
Choose $\eta\in C_c^\infty((1/2,2))$ with $\sum_{\lambda\in2^{\mathbb Z}}\eta(\p/\lambda)=1$ for $\p>0$.
For $0<\lambda\le2C_m$, set
\begin{equation*}
I^\pm_\lambda(t;r_1,r_2)
:=
\int_0^\infty
e^{\pm it\sqrt{\p^2+m^2}}
J_\nu(\p r_1)J_\mu(\p r_2)\chi(\p)\eta(\p/\lambda)\,\p\,d\p .
\end{equation*}
After the change of variables $\p=\lambda s$, we get
\begin{equation*}
I^\pm_\lambda
=
\lambda^2 e^{\pm imt}
\int_0^\infty
e^{\pm i\lambda^2t\,\Omega_\lambda(s)}
b_\lambda(s)
J_\nu(\lambda r_1s)J_\mu(\lambda r_2s)s\,ds,
\end{equation*}
where
\begin{equation*}
        \Omega_\lambda(s)
        :=
        \lambda^{-2}\bigl(\sqrt{m^2+\lambda^2s^2}-m\bigr),
\end{equation*}
and $b_\lambda(s)=\chi(\lambda s)\eta(s)$ is supported in a fixed compact subset of
$(0,\infty)$, with uniform $C^N$-bounds. On this support,
\begin{equation*}
        |\Omega_\lambda'(s)|\ge c>0,
        \qquad
        |\partial_s^N\Omega_\lambda(s)|\le C_N,\qquad N\ge2,
\end{equation*}
uniformly for $0<\lambda\le2C_m$. Therefore Lemma~\ref{lem:fixed-shell-bessel-massive},
applied with
\begin{equation*}
        T=\lambda^2|t|,
        \qquad
        R_1=\lambda r_1,
        \qquad
        R_2=\lambda r_2,
\end{equation*}
gives
\begin{equation*}
|I^\pm_\lambda(t;r_1,r_2)|
\le
C\lambda^2(1+\lambda^2|t|)^{-1/2}
(1+(\lambda r_1)^{-\nu_-})(1+(\lambda r_2)^{-\mu_-}).
\end{equation*}
Since, for $\sigma\ge0$,
\begin{equation*}
        1+(\lambda r_i)^{-\sigma}
        \le C\lambda^{-\sigma}(1+r_i^{-\sigma}),
        \qquad 0<\lambda\le2C_m,
\end{equation*}
we obtain
\begin{equation*}
\sum_{\lambda\le2C_m}|I^\pm_\lambda(t;r_1,r_2)|
\le
C(1+r_1^{-\nu_-})(1+r_2^{-\mu_-})
\sum_{\lambda\le2C_m}
\lambda^{2-\nu_- -\mu_-}(1+\lambda^2|t|)^{-1/2}.
\end{equation*}
Writing $a=\nu_-+\mu_-<1$, the last dyadic sum satisfies
\begin{equation*}
\sum_{\lambda\le2C_m}
\lambda^{2-a}(1+\lambda^2|t|)^{-1/2}
\le C|t|^{-1/2}.
\end{equation*}
Indeed, the part $\lambda\le |t|^{-1/2}$ is bounded by
$C|t|^{-(2-a)/2}\le C|t|^{-1/2}$, while the part
$\lambda>|t|^{-1/2}$ is bounded by
\begin{equation*}
C|t|^{-1/2}\sum_{\lambda\le2C_m}\lambda^{1-a}
\le C|t|^{-1/2},
\end{equation*}
because $a<1$. Hence
\begin{equation*}
|I^\pm_{\nu,\mu}(t;r_1,r_2)|
\le
C|t|^{-1/2}
(1+r_1^{-\nu_-})(1+r_2^{-\mu_-})
\end{equation*}
whenever $\nu_-+\mu_-<1$.

It remains to treat the borderline case
\begin{equation*}
        \nu=\mu=-\frac12 .
\end{equation*}
Here we use the explicit formula
\begin{equation*}
        J_{-1/2}(z)=\sqrt{\frac2{\pi z}}\cos z .
\end{equation*}
Thus
\begin{equation*}
J_{-1/2}(\p r_1)J_{-1/2}(\p r_2)\p
=
\frac2\pi (r_1r_2)^{-1/2}\cos(\p r_1)\cos(\p r_2),
\end{equation*}
and $I^\pm_{-1/2,-1/2}$ is a finite sum of terms of the form
\begin{equation*}
(r_1r_2)^{-1/2}
\int_0^\infty
e^{i\Phi_{\varepsilon_1,\varepsilon_2}(\p)}
\chi(\p)\,d\p,
\end{equation*}
where
\begin{equation*}
        \Phi_{\varepsilon_1,\varepsilon_2}(\p)
        =
        \pm t\sqrt{\p^2+m^2}
        +\varepsilon_1\p r_1+\varepsilon_2\p r_2,
        \qquad
        \varepsilon_1,\varepsilon_2\in\{-1,1\}.
\end{equation*}
On $\operatorname{supp}\chi$,
\begin{equation*}
        |\Phi_{\varepsilon_1,\varepsilon_2}''(\p)|
        =
        |t|\frac{m^2}{(\p^2+m^2)^{3/2}}
        \ge c_m|t|.
\end{equation*}
Van der Corput's lemma gives
\begin{equation*}
        |I^\pm_{-1/2,-1/2}(t;r_1,r_2)|
        \le
        C|t|^{-1/2}(r_1r_2)^{-1/2}.
\end{equation*}
Since
\begin{equation*}
        (r_1r_2)^{-1/2}
        \le
        (1+r_1^{-1/2})(1+r_2^{-1/2}),
\end{equation*}
this is the desired estimate in the endpoint case.

Combining the cases $|t|\le1$ and $|t|\ge1$ concludes the proof.
\end{proof}

Before getting to the proof of estimate \eqref{est:dis0}, let us fix the setup, returning to the four kernels \eqref{eq:Fj-massive}--\eqref{eq:Stildej-massive}.  After
changing variable from $\rho$ to $\p=\sqrt{\rho^{2}-m^{2}}$, we have
\begin{equation*}
\rho=\sqrt{\p^{2}+m^{2}},
\qquad
\frac{d\rho}{d\p}=\frac{\p}{\rho}.
\end{equation*}
Hence, writing
\begin{equation*}
a^F_{j,\pm}(\p)
:=
\frac{\sqrt{\p^2+m^2}\mp m}{4\pi\sqrt{\p^2+m^2}}
\varphi_j(\sqrt{\p^2+m^2}),
\end{equation*}

\begin{equation*}
a^E_{j,\pm}(\p)
:=
\frac{\sqrt{\p^2+m^2}\pm m}{4\pi\sqrt{\p^2+m^2}}
\varphi_j(\sqrt{\p^2+m^2}).
\end{equation*}
\begin{equation*}
a_{j,\pm}^{S}(\p):=\pm\frac{\p}{4\pi i\sqrt{\p^{2}+m^{2}}}
\varphi_{j}(\sqrt{\p^{2}+m^{2}}),
\end{equation*}
we use the following amplitude notation at every energy:
\begin{equation*}
I_{\nu,\mu}^{\pm}[a](t;r_1,r_2)
:=\int_0^\infty e^{\pm it\sqrt{\p^2+m^2}}
J_\nu(\p r_1)J_\mu(\p r_2)a(\p)\p\,d\p.
\end{equation*}
The kernels can then be written as
\begin{equation}\label{eq:F-model-massive}
F_{\alpha,j}^{m}
=
I^{+}_{-\alpha,-\alpha}[a_{j,+}^{F}]
+
I^{-}_{-\alpha,-\alpha}[a_{j,-}^{F}],
\end{equation}
\begin{equation}\label{eq:E-model-massive}
E_{\alpha,j}^{m}
=
I^{+}_{1-\alpha,1-\alpha}[a_{j,+}^{E}]
+
I^{-}_{1-\alpha,1-\alpha}[a_{j,-}^{E}],
\end{equation}
\begin{equation}\label{eq:S-model-massive}
S_{\alpha,j}^{m}
=
I^{+}_{-\alpha,1-\alpha}[a_{j,+}^{S}]
+
I^{-}_{-\alpha,1-\alpha}[a_{j,-}^{S}],
\end{equation}
\begin{equation}\label{eq:Stilde-model-massive}
\widetilde S_{\alpha,j}^{m}
=
I^{+}_{1-\alpha,-\alpha}[a_{j,+}^{S}]
+
I^{-}_{1-\alpha,-\alpha}[a_{j,-}^{S}].
\end{equation}

Then we have the following

\begin{lemma}\label{lem:kernelbounds-massive-sing}
There exists a constant $C_{\alpha,m}$ such that for all $j\in\mathbb{Z}$,
$t\in\mathbb{R}$ and $r_{1},r_{2}>0$ one has
\begin{equation}\label{eq:E-bound-massive}
|E_{\alpha,j}^{m}(t;r_{1},r_{2})|
\le
C_{\alpha,m}2^{2j}(1+2^{j}|t|)^{-1/2},
\end{equation}
\begin{equation}\label{eq:F-bound-massive}
|F_{\alpha,j}^{m}(t;r_{1},r_{2})|
\le
C_{\alpha,m}2^{2j}(1+2^{j}|t|)^{-1/2}
\big(1+(2^{j}r_{1})^{-\alpha}\big)
\big(1+(2^{j}r_{2})^{-\alpha}\big),
\end{equation}
\begin{equation}\label{eq:S-bound-massive}
|S_{\alpha,j}^{m}(t;r_{1},r_{2})|
\le
C_{\alpha,m}2^{2j}(1+2^{j}|t|)^{-1/2}
\big(1+(2^{j}r_{1})^{-\alpha}\big),
\end{equation}
\begin{equation}\label{eq:Stilde-bound-massive}
|\widetilde S_{\alpha,j}^{m}(t;r_{1},r_{2})|
\le
C_{\alpha,m}2^{2j}(1+2^{j}|t|)^{-1/2}
\big(1+(2^{j}r_{2})^{-\alpha}\big).
\end{equation}
\end{lemma}

\begin{proof}
If $2^{j}<m/2$, then all four kernels vanish identically, since $\varphi_{j}(\rho)=0$ on
$[m,\infty)$. If $m/2\le2^{j}\le4m$, there are only finitely many such $j$, and the
estimates follow from Lemma~\ref{lem:threshold-massive-sing} and from the
representations \eqref{eq:F-model-massive}--\eqref{eq:Stilde-model-massive}. Finally, if $2^j\ge4m$, we set $\p=2^j\lambda$ and
$\tau_j=m2^{-j}\le1/4$. The rescaled phase is
\begin{equation*}
\omega_{\tau_j}(\lambda)=\sqrt{\lambda^2+\tau_j^2},
\end{equation*}
which satisfies the assumptions of
Lemma~\ref{lem:fixed-shell-bessel-massive} uniformly in $j$.
Rescaling back gives the required factor
\begin{equation*}
2^{2j}(1+2^j|t|)^{-1/2}.
\end{equation*}\end{proof}
We are finally ready to prove estimate \eqref{est:dis0}.

\begin{proof}
Let
\begin{equation*}
        f_j=\widetilde\varphi_j(|\mathcal{D}_A|)P_0f
        =
        \binom{\phi_j(r)}{e^{i\theta}\psi_j(r)}
\end{equation*}
and set
\begin{equation*}
        w_j(r):=(1+(2^jr)^{-\alpha})^{-1}.
\end{equation*}
Then
\begin{equation*}
        W_j f_j
        =
        \binom{w_j(r)\phi_j(r)}{e^{i\theta}\psi_j(r)}.
\end{equation*}
By the kernel representation \eqref{eq:singular-localized-kernel-massive}, applied to
the datum $W_jf_j$, one has
\begin{equation*}
\varphi_j(|D_A|)e^{itD_A}P_0W_jf_j
=
\binom{
T^F_j(w_j\phi_j)-T^S_j\psi_j
}{
e^{i\theta}\bigl(T^{\widetilde S}_j(w_j\phi_j)+T^E_j\psi_j\bigr)
},
\end{equation*}
where
\begin{equation*}
(T^F_j g)(r):=2\pi\int_0^\infty F^m_{\alpha,j}(t;r,r_2)g(r_2)\,r_2dr_2,
\end{equation*}
\begin{equation*}
(T^E_j g)(r):=2\pi\int_0^\infty E^m_{\alpha,j}(t;r,r_2)g(r_2)\,r_2dr_2,
\end{equation*}

\begin{equation*}
\begin{aligned}
(T^S_j g)(r)&:=2\pi\int_0^\infty S^m_{\alpha,j}(t;r,r_2)g(r_2)\,r_2dr_2,\\
(T^{\widetilde S}_j g)(r)&:=2\pi\int_0^\infty \widetilde S^m_{\alpha,j}(t;r,r_2)g(r_2)\,r_2dr_2.
\end{aligned}
\end{equation*}
By Lemma~\ref{lem:kernelbounds-massive-sing},
\begin{equation*}
\begin{aligned}
\|w_jT^F_j(w_j\phi_j)\|_{L^\infty(\mathbb R^2)}
&\lesssim
2^{2j}(1+2^j|t|)^{-1/2}
\|(1+(2^jr)^{-\alpha})w_j\phi_j\|_{L^1(\mathbb R^2)}
\\
&=
2^{2j}(1+2^j|t|)^{-1/2}
\|\phi_j\|_{L^1(\mathbb R^2)}.
\end{aligned}
\end{equation*}
Similarly,
\begin{equation*}
\begin{aligned}
\|T^{\widetilde S}_j(w_j\phi_j)\|_{L^\infty(\mathbb R^2)}
&\lesssim
2^{2j}(1+2^j|t|)^{-1/2}
\|(1+(2^jr)^{-\alpha})w_j\phi_j\|_{L^1(\mathbb R^2)}
\\
&=
2^{2j}(1+2^j|t|)^{-1/2}
\|\phi_j\|_{L^1(\mathbb R^2)}.
\end{aligned}
\end{equation*}
For the remaining two terms, again by Lemma~\ref{lem:kernelbounds-massive-sing},
\begin{equation*}
\|w_jT^S_j\psi_j\|_{L^\infty(\mathbb R^2)}
+
\|T^E_j\psi_j\|_{L^\infty(\mathbb R^2)}
\lesssim
2^{2j}(1+2^j|t|)^{-1/2}
\|\psi_j\|_{L^1(\mathbb R^2)}.
\end{equation*}
Combining the four estimates gives
\begin{equation*}
\left\|
W_j e^{it\mathcal{D}_A}\varphi_j(|\mathcal{D}_A|)P_0W_j f_j
\right\|_{L^\infty(\mathbb R^2)}
\lesssim
2^{2j}(1+2^j|t|)^{-1/2}
\|f_j\|_{L^1(\mathbb R^2)}.
\end{equation*}
Since $f_j=\widetilde\varphi_j(|\mathcal{D}_A|)P_0f$, this is precisely
\eqref{est:dis0}.
\end{proof}

\subsection{Unweighted estimate in the singular sector: proof of \eqref{est:disq}}\label{subsec:unweighted-dispersive}

We now prove the localized unweighted estimate in the singular channel.
The argument is the massive analogue of Subsection~5.3 from \cite{cacdanyinzha}.
As in the erratum, besides the diagonal singular kernel one has to control the two
mixed off diagonal terms; the key point is that in the mixed kernels only one Bessel
factor has negative order, so they are no worse than the diagonal singular one.

After the change of variable
\begin{equation*}
       \p=\sqrt{\rho^2-m^2},
\end{equation*}
the singular kernels are finite linear combinations of oscillatory Bessel kernels with
phase
\begin{equation*}
        e^{\pm it\sqrt{\p^2+m^2}}.
\end{equation*}
The only reusable model estimate we need is the high frequency dyadic one below; the
finite threshold range $m/2\le2^j\le4m$ will be treated directly in the proof of
\eqref{est:disq}.

Let us start with the following
\begin{lemma}\label{cor:model-Lq-dyadic-massive}
Let $-\frac12\le \nu<0$ and assume either $\mu=\nu$ or $\mu\ge0$.
Let $a_j\in C^\infty((0,\infty))$ satisfy
\begin{equation*}
\supp a_j\subset[2^{j-2},2^{j+2}],
\qquad
|\partial_\p^N a_j(\p)|\le C_N\,2^{-Nj},
\qquad N\ge0,
\end{equation*}
and assume $2^j\ge4m$. Define
\begin{equation*}
I_{\nu,\mu,j}^{\pm}(t;r_1,r_2)
:=
\int_0^\infty
e^{\pm it\sqrt{\p^2+m^2}}
a_j(\p)\,J_\nu(\p r_1)J_\mu(\p r_2)\,\p\,d\p
\end{equation*}
and
\begin{equation*}
(T_{\nu,\mu,j}^{\pm}g)(t,r_1)
:=
\int_0^\infty
I_{\nu,\mu,j}^{\pm}(t;r_1,r_2)\,g(r_2)\,r_2\,dr_2.
\end{equation*}
Then for every $2\le q<-\frac2\nu$ one has
\begin{equation}\label{eq:model-Lq-dyadic-massive}
\|T_{\nu,\mu,j}^{\pm}g(t,\cdot)\|_{L^q((0,\infty),rdr)}
\le
C\,2^{2j(1-\frac2q)}(1+2^j|t|)^{-\frac12(1-\frac2q)}
\|g\|_{L^{q'}((0,\infty),rdr)},
\end{equation}
and the same estimate holds with $(\nu,\mu)$ replaced by $(\mu,\nu)$.
\end{lemma}

\begin{proof}
Set $N=2^j$, $\tau=m/N\leq1/4$, $s=Nt$, and
$b_j(\lambda)=a_j(N\lambda)$. The rescaled kernels are
\begin{equation*}
K_j^\pm(s;R,S)
:=
\int_0^\infty
 e^{\pm is\omega_\tau(\lambda)}b_j(\lambda)
 J_\nu(R\lambda)J_\mu(S\lambda)\lambda\,d\lambda,
\qquad
\omega_\tau(\lambda)=\sqrt{\lambda^2+\tau^2}.
\end{equation*}
Let $\widetilde T_j^\pm(s)$ denote the corresponding integral operator
with measure $S\,dS$. Throughout the proof, spatial $L^q$ norms are
taken with respect to the radial measure.

The symbols $b_j$ are supported in $[1/4,4]$ and have uniformly
bounded derivatives. On this interval, $\omega_\tau'$ is bounded
away from zero and all its derivatives are uniformly bounded.
Lemma~\ref{lem:fixed-shell-bessel-massive} therefore gives
\begin{equation}\label{eq:high-Lq-pointwise}
|K_j^\pm(s;R,S)|
\lesssim
(1+|s|)^{-1/2}
(1+R^{-\nu_-})(1+S^{-\mu_-}),
\end{equation}
where $\nu_-:=\max\{0,-\nu\}$ and
$\mu_-:=\max\{0,-\mu\}$. Also, by Hankel--Plancherel,
\begin{equation*}
\widetilde T_j^\pm(s)h
=
\mathcal H_\nu
\bigl[e^{\pm is\omega_\tau}b_j\mathcal H_\mu h\bigr],
\qquad
\|\widetilde T_j^\pm(s)\|_{L^2\to L^2}\lesssim1.
\end{equation*}
All constants are independent of $j$ and $s$.

We use the spatial decomposition of
\cite[Lemma~5.4]{cacdanyinzha}, setting
$\chi_0=\mathbf 1_{(0,1)}$ and $\chi_1=1-\chi_0$.
Since $q<-2/\nu$ and either $\mu=\nu$ or $\mu\geq0$, both
$\chi_0(R)(1+R^{-\nu_-})$ and
$\chi_0(R)(1+R^{-\mu_-})$ belong to $L^q(R\,dR)$.
Thus \eqref{eq:high-Lq-pointwise} and H\"older's inequality give
\begin{equation*}
\|\chi_0\widetilde T_j^\pm(s)\chi_0\|_{L^{q'}\to L^q}
\lesssim (1+|s|)^{-1/2}.
\end{equation*}
For the large--large piece, \eqref{eq:high-Lq-pointwise} gives the
$L^1\to L^\infty$ bound with decay $(1+|s|)^{-1/2}$.
Interpolating with its uniform $L^2\to L^2$ bound yields
\begin{equation*}
\|\chi_1\widetilde T_j^\pm(s)\chi_1\|_{L^{q'}\to L^q}
\lesssim (1+|s|)^{-\frac12(1-\frac2q)}.
\end{equation*}

For the mixed piece $0<R<1\leq S$, we apply Hankel--Plancherel
to the kernel as a function of $S$. The compact frequency support
and the small-argument Bessel bound give
\begin{equation*}
\begin{aligned}
\|K_j^\pm(s;R,\cdot)\|_{L^2((0,\infty),S\,dS)}
&=
\|b_j(\lambda)J_\nu(R\lambda)\|_{L^2(\lambda\,d\lambda)}
\lesssim R^\nu.
\end{aligned}
\end{equation*}
On the other hand, \eqref{eq:high-Lq-pointwise} implies
\begin{equation*}
\|K_j^\pm(s;R,\cdot)\|_{L^\infty(1,\infty)}
\lesssim R^\nu(1+|s|)^{-1/2}.
\end{equation*}
Interpolating these two norms, for $2\leq q<-2/\nu$, gives
\begin{equation*}
\|K_j^\pm(s;R,\cdot)\|_{L^q((1,\infty),S\,dS)}
\lesssim
R^\nu(1+|s|)^{-\frac12(1-\frac2q)}.
\end{equation*}
Consequently, H\"older's inequality and
$\int_0^1R^{\nu q+1}\,dR<\infty$ yield
\begin{equation*}
\|\chi_0\widetilde T_j^\pm(s)\chi_1h\|_{L^q}
\lesssim
(1+|s|)^{-\frac12(1-\frac2q)}\|h\|_{L^{q'}}.
\end{equation*}
The opposite mixed piece is treated by applying the same argument
to its adjoint, whose orders are $(\mu,\nu)$, with conjugate
amplitude and opposite phase sign. If $\mu\geq0$, the small-radius
factors in the two kernel-row bounds are simply bounded, so no
additional restriction on $q$ occurs.

Summing the four pieces, we obtain
\begin{equation*}
\|\widetilde T_j^\pm(s)h\|_{L^q}
\lesssim
(1+|s|)^{-\frac12(1-\frac2q)}\|h\|_{L^{q'}}.
\end{equation*}
Finally, put $g_N(S)=g(S/N)$. The change of variables
$\p=N\lambda$ gives
\begin{equation*}
(T_{\nu,\mu,j}^\pm g)(t,R/N)
=
\bigl(\widetilde T_j^\pm(Nt)g_N\bigr)(R),
\qquad
\|g_N\|_{L^{q'}}=N^{2/q'}\|g\|_{L^{q'}}.
\end{equation*}
It follows that
\begin{equation*}
\|T_{\nu,\mu,j}^\pm g(t,\cdot)\|_{L^q}
\lesssim
N^{2/q'-2/q}(1+N|t|)^{-\frac12(1-\frac2q)}
\|g\|_{L^{q'}},
\end{equation*}
which is \eqref{eq:model-Lq-dyadic-massive}, since
$2/q'-2/q=2(1-2/q)$.
The estimate with the orders interchanged follows by the same
adjoint argument.
\end{proof}
We are now in a position to prove estimate \eqref{est:disq}.

\begin{proof}
Write
\begin{equation*}
f_j:=\widetilde\varphi_j(|\mathcal{D}_A|)P_0f
=
\begin{pmatrix}
g_0(r)\\
e^{i\theta}g_1(r)
\end{pmatrix}.
\end{equation*}
By the explicit form of the localized singular propagator established in
Subsection~\ref{subsec:weighted-dispersive},
\begin{equation*}
\varphi_j(|\mathcal{D}_A|)e^{it\mathcal{D}_A}P_0f
=
\begin{pmatrix}
T_{F,j}^m g_0 - T_{S,j}^{12,m}g_1\\[1mm]
e^{i\theta}\bigl(T_{S,j}^{21,m}g_0 + T_{E,j}^m g_1\bigr)
\end{pmatrix},
\end{equation*}
where
\begin{equation*}
\begin{aligned}
(T_{F,j}^m h)(r)&=2\pi\int_0^\infty F_{\alpha,j}^m(t;r,r_2)h(r_2)\,r_2dr_2,\\
(T_{E,j}^m h)(r)&=2\pi\int_0^\infty E_{\alpha,j}^m(t;r,r_2)h(r_2)\,r_2dr_2,
\end{aligned}
\end{equation*}
and similarly for $T_{S,j}^{12,m}$ and $T_{S,j}^{21,m}$.

We first treat the regular diagonal term $T_{E,j}^m$. Since $E_{\alpha,j}^m$
contains only positive-order Bessel factors, the $L^1\to L^\infty$ estimate proved in
the kernel bound \eqref{eq:E-bound-massive} gives
\begin{equation*}
\|T_{E,j}^m h\|_{L^\infty((0,\infty),rdr)}
\le
C\,2^{2j}(1+2^j|t|)^{-1/2}\|h\|_{L^1((0,\infty),rdr)}.
\end{equation*}
On the other hand, $T_{E,j}^m$ is a scalar spectral multiplier in the Hankel variable
of order $1-\alpha$, hence Hankel--Plancherel yields
\begin{equation*}
\|T_{E,j}^m h\|_{L^2((0,\infty),rdr)}\le C\|h\|_{L^2((0,\infty),rdr)}.
\end{equation*}
Interpolating these two bounds gives
\begin{equation}\label{eq:E-Lq-massive}
\|T_{E,j}^m h\|_{L^q((0,\infty),rdr)}
\le
C\,2^{2j(1-\frac2q)}(1+2^j|t|)^{-\frac12(1-\frac2q)}
\|h\|_{L^{q'}((0,\infty),rdr)}.
\end{equation}
We next treat the singular diagonal term $T_{F,j}^m$ and the mixed terms
$T_{S,j}^{12,m}$, $T_{S,j}^{21,m}$.

\smallskip
\noindent
\emph{1). The low-frequency range $2^j<m/2$.}
In this case all three kernels $F_{\alpha,j}^m$, $S_{\alpha,j}^{12,m}$,
$S_{\alpha,j}^{21,m}$ vanish identically, since $\varphi_j(\rho)=0$ on $[m,\infty)$.

\smallskip
\noindent
\emph{2). The threshold range $m/2\le 2^j\le 4m$}.
There are only finitely many such $j$.
All estimates in this range are therefore performed at a fixed physical scale, and the
constants are uniform in $j$. The only genuinely new point with respect to the
massless proof is that, after the change of variable
\begin{equation*}
        \p=\sqrt{\rho^2-m^2},
\end{equation*}
the oscillatory phase becomes
\begin{equation*}
        e^{\pm it\sqrt{\p^2+m^2}}.
\end{equation*}
Near $\p=0$ this phase is quadratic. We can nevertheless use the same
spatial decomposition as at high frequency, without an additional split
depending on $\p r_2$. Indeed, the kernels are finite sums of operators
\begin{equation*}
T=\mathcal H_\nu
M_{e^{\pm it\sqrt{\p^2+m^2}}a_j(\p)}\mathcal H_\mu,
\end{equation*}
where $a_j$ is smooth and supported in a fixed interval $[0,C_m]$, and
$(\nu,\mu)$ is one of $(-\alpha,-\alpha)$,
$(-\alpha,1-\alpha)$ and $(1-\alpha,-\alpha)$.
Write $K(t;r_1,r_2)$ for its kernel and set
$\chi_0=\mathbf1_{(0,1)}$, $\chi_1=1-\chi_0$.
Lemma~\ref{lem:threshold-massive-sing} gives

\begin{equation*}
|K(t;r_1,r_2)|\lesssim(1+|t|)^{-1/2}
(1+r_1^{-\nu_-})(1+r_2^{-\mu_-}).
\end{equation*}

The local--local part is bounded from $L^{q'}$ to $L^q$ by H\"older's
inequality, since $q<2/\alpha$. The large--large part has this
$L^1\to L^\infty$ bound and a uniform $L^2\to L^2$ bound by
Hankel--Plancherel. Interpolation gives the required decay
$(1+|t|)^{-\frac12(1-2/q)}$.

For a mixed part, fix $0<r_1<1$. Hankel--Plancherel in $r_2$ gives
\begin{equation*}
\begin{aligned}
\|K(t;r_1,\cdot)\|_{L^2(r_2dr_2)}^2
&=\int_0^{C_m}|a_j(\p)J_\nu(\p r_1)|^2\p\,d\p\\
&\lesssim r_1^{2\nu}\int_0^{C_m}\p^{2\nu+1}d\p
\lesssim r_1^{2\nu},
\end{aligned}
\end{equation*}
because $\nu\ge-1/2>-1$. On $r_2\ge1$, the pointwise threshold bound
also gives
\begin{equation*}
\|K(t;r_1,\cdot)\|_{L^\infty(1,\infty)}
\lesssim (1+|t|)^{-1/2}(1+r_1^{-\nu_-}).
\end{equation*}
Interpolating these two row norms and using $r_1^\nu\lesssim
1+r_1^{-\nu_-}$ on $(0,1)$ yields
\begin{equation*}
\|K(t;r_1,\cdot)\|_{L^q((1,\infty),r_2dr_2)}
\lesssim (1+|t|)^{-\frac12(1-2/q)}(1+r_1^{-\nu_-}).
\end{equation*}
Integrating in $r_1$, followed by H\"older's inequality in $r_2$,
proves the mixed bound. The opposite mixed part follows by taking
the adjoint and interchanging $\nu$ and $\mu$.

Summing the four spatial pieces, we obtain
\begin{equation}\label{eq:threshold-Lq-massive}
\|T_{F,j}^m h\|_{L^q}
+
\|T_{S,j}^{12,m} h\|_{L^q}
+
\|T_{S,j}^{21,m} h\|_{L^q}
\le
C(1+|t|)^{-\frac12(1-\frac2q)}
\|h\|_{L^{q'}},
\end{equation}
with $C$ uniform in the threshold range. Since here $2^j\simeq m$, estimate
\eqref{eq:threshold-Lq-massive} is exactly equivalent to the bound required in
\eqref{est:disq}.

\smallskip
\noindent
\emph{3). The high frequency range $2^j\ge4m$.}
After the same change of variable, the kernels $F_{\alpha,j}^m$,
$S_{\alpha,j}^{12,m}$ and $S_{\alpha,j}^{21,m}$ are finite linear combinations of
dyadic model operators satisfying the assumptions of
Lemma~\ref{cor:model-Lq-dyadic-massive}, with the same
three order pairs listed above and with amplitudes satisfying the assumptions of
Lemma~\ref{cor:model-Lq-dyadic-massive}. Hence
\begin{equation}\label{eq:high-Lq-massive}
\|T_{F,j}^m h\|_{L^q}
+
\|T_{S,j}^{12,m} h\|_{L^q}
+
\|T_{S,j}^{21,m} h\|_{L^q}
\le
C\,2^{2j(1-\frac2q)}(1+2^j|t|)^{-\frac12(1-\frac2q)}
\|h\|_{L^{q'}}.
\end{equation}

Combining \eqref{eq:E-Lq-massive}, \eqref{eq:threshold-Lq-massive} and
\eqref{eq:high-Lq-massive} with the decomposition of the localized singular flow proves
\eqref{est:disq}.
\end{proof}

\section{Strichartz estimates}\label{sec:strichartz}

\subsection{Strichartz for the regular component: proof of Theorem \ref{th:strichartz-regular-massive} }

\begin{proof}
The case $(p,q)=(\infty,2)$ follows immediately from the spectral theorem.

We therefore assume $\eta\in[0,1]$ and
\begin{equation*}
(p,q)\in(2,\infty)^2,
\qquad
\frac2p+\frac{1+\eta}{q}\leq\frac{1+\eta}{2},
\qquad
s=(2+\eta)\left(\frac12-\frac1q\right)-\frac1p.
\end{equation*}

We use the scalar Klein--Gordon Strichartz estimate associated with the
Friedrichs Aharonov--Bohm Hamiltonian. More precisely, by
\cite[Definition~1.2 and Theorem~4.1]{GaoYinZhangZheng}, for either choice of sign and for
arbitrary two-component data $g$,
\begin{equation}\label{eq:scalar-KG-Strichartz}
\big\|
e^{\pm it\sqrt{\mathscr L_m}}g
\big\|_{L_t^p(\R;[L_x^q(\R^2)]^2)}
\lesssim
\big\|(I+\mathscr L_m)^{s/2}g\big\|_{[L^2(\R^2)]^2}.
\end{equation}
The result in \cite{GaoYinZhangZheng} is stated after normalizing the
mass; the estimate for an arbitrary fixed $m>0$ follows by dilation,
using the homogeneity of the Aharonov--Bohm Hamiltonian.

We now apply \eqref{eq:scalar-KG-Strichartz} to
$g=\Pi_\perp^\pm f$.  By \eqref{eq:regular-splitting-massive},
\begin{equation*}
\begin{aligned}
\|e^{it\mathcal D_A}P_\perp f\|_{L_t^p[L_x^q]^2}
&\leq
\|e^{it\sqrt{\mathscr L_m}}\Pi_\perp^+f\|_{L_t^p[L_x^q]^2}
+
\|e^{-it\sqrt{\mathscr L_m}}\Pi_\perp^-f\|_{L_t^p[L_x^q]^2}
\\
&\lesssim
\big\|(I+\mathscr L_m)^{s/2}\Pi_\perp^+f\big\|_{[L^2]^2}
+
\big\|(I+\mathscr L_m)^{s/2}\Pi_\perp^-f\big\|_{[L^2]^2}.
\end{aligned}
\end{equation*}
On $\mathcal H_\perp$ one has
\begin{equation*}
\mathscr L_m=D_\perp^2.
\end{equation*}
Hence
\begin{equation*}
\big\|(I+\mathscr L_m)^{s/2}g\big\|_{[L^2]^2}
=
\|\langle\mathcal D_A\rangle^s g\|_{[L^2]^2}
=
\|g\|_{H_A^s},
\qquad g\in\mathcal H_\perp.
\end{equation*}
Moreover, $\Pi_\perp^\pm$ are mutually orthogonal spectral projections of
$D_\perp$ and commute with $\langle D_\perp\rangle^s$. Therefore
\begin{equation*}
\|\Pi_\perp^+f\|_{H_A^s}^2
+
\|\Pi_\perp^-f\|_{H_A^s}^2
=
\|P_\perp f\|_{H_A^s}^2.
\end{equation*}
We conclude
\begin{equation*}
\|e^{it\mathcal D_A}P_\perp f\|_{L_t^p(\R;[L_x^q(\R^2)]^2)}
\lesssim
\|P_\perp f\|_{H_A^s(\R^2;\C^2)},
\end{equation*}
which proves the theorem.
\end{proof}

\begin{remark}
Notice that this proof is clearly independent of the choice of the self-adjoint extension of $\mathcal{D}_A$, as indeed that only affects the critical sector $k=0$, which does not play any role in Theorem \ref{th:strichartz-regular-massive}.
\end{remark}

\subsection{Unweighted Strichartz estimates for the singular component}
We use a square function inequality and a localized spacetime estimate.
The square function bound is proved directly below, with a constant
uniform in all finite randomized sums. The localized estimate uses the
linear time scale at high momentum and the quadratic time scale near
the massive thresholds.

By reflection it suffices to work with $0<\alpha<1/2$ and the selected
extension $\gamma=3\pi/4$; we write $\mathcal D_A=\mathcal D_{A,\gamma}$.
In particular, $q(\alpha)=2/\alpha>4$. The other flux range follows by
replacing $\alpha$ with $1-\alpha$ and interchanging the components.
Note that the following result and its proof work also in the
massless case $m=0$
(just take $\omega(\xi)=\xi$ in the argument).

\begin{proposition}[Square function inequality on the singular channel]
\label{prop:square function-singular-massive}
Let $\alpha\in(0,\frac12)$ and let $0\le\varphi\le1$ be a smooth function
supported in $[\frac12,2]$ satisfying
\begin{equation*}
\sum_{j\in\Z}\varphi(2^{-j}\lambda)=1,
\qquad \lambda>0.
\end{equation*}
With $q(\alpha)$ defined in \eqref{q-alp}, set
\begin{equation*}
q'(\alpha):=\frac{q(\alpha)}{q(\alpha)-1}
=\frac{2}{2-\alpha}.
\end{equation*}
Then for every
\begin{equation*}
q'(\alpha)<r<q(\alpha)
\end{equation*}
there exist constants $c_r,C_r>0$ such that
\begin{equation}\label{eq:square function-singular-massive}
c_r\|P_0f\|_{[L^r(\R^2)]^2}
\le
\Big\|
\Big(\sum_{j\in\Z}|\varphi_j(|\mathcal D_A|)P_0f|^2\Big)^{1/2}
\Big\|_{[L^r(\R^2)]^2}
\le
C_r\|P_0f\|_{[L^r(\R^2)]^2}.
\end{equation}
\end{proposition}

\begin{proof}
Put $\omega(\xi)=\sqrt{m^2+\xi^2}$. For
$P_0f(x,\theta)=(f_1(x),e^{i\theta}f_2(x))^T$, the spectral
representation gives, for every bounded Borel function $F$,
\begin{equation}\label{eq:massive-functional-calculus-zero-mode}
F(|\mathcal D_A|)P_0f
=
\binom{
\mathcal H_{-\alpha}
\bigl[F(\omega(\p))\mathcal H_{-\alpha}f_1\bigr]
}{
e^{i\theta}
\mathcal H_{1-\alpha}
\bigl[F(\omega(\p))\mathcal H_{1-\alpha}f_2\bigr]
}.
\end{equation}
Thus it suffices to prove the scalar square function estimates for the
orders $\nu=-\alpha$ and $\nu=1-\alpha$.

We first establish a multiplier bound. For $-1/2<\nu\le1$, put
\begin{equation*}
 d=2\nu+2,\qquad d\mu_\nu(x)=x^{2\nu+1}\,dx,\qquad
 \widetilde J_\nu(z)=z^{-\nu}J_\nu(z).
\end{equation*}
The modified Hankel transform
\begin{equation*}
 h_\nu g(\xi)=\int_0^\infty\widetilde J_\nu(x\xi)g(x)\,d\mu_\nu(x)
\end{equation*}
is an involutive isometry on $L^2(d\mu_\nu)$. We claim that
\begin{equation}\label{eq:claim}
 \|h_\nu a h_\nu g\|_{L^r(x^\gamma d\mu_\nu)}
 \lesssim_{\nu,r,\gamma} A_4(a)\|g\|_{L^r(x^\gamma d\mu_\nu)},
 \qquad -d<\gamma<d(r-1),
\end{equation}
where $1<r<\infty$, $a\in C^4(0,\infty)$ and
\begin{equation*}
 A_4(a)=\max_{0\le k\le4}
 \|(\xi\partial_\xi)^k a\|_\infty<\infty.
\end{equation*}

Choose $\eta\in C_c^\infty((1/2,2))$ such that $\eta\ge0$,
$\sum_{N\in2^{\mathbb Z}}\eta(\xi/N)=1$ for $\xi>0$, and set
\begin{equation*}
 a_N(\xi)=a(\xi)\eta(\xi/N),\qquad
 k_N=h_\nu a_N,\qquad v_N(s)=a(Ns)\eta(s).
\end{equation*}
The $C^4$ norms of $v_N$ are bounded by $C A_4(a)$,
since $s$ stays in a fixed annulus. The bounds in
Lemma~\ref{lem:Bessel-massive-sing}, together with
$\widetilde J_\nu'(z)=-z^{-\nu}J_{\nu+1}(z)$, give
\begin{equation*}
 |\widetilde J_\nu(z)|+|\widetilde J_\nu'(z)|
 \lesssim_\nu (1+z)^{-(d-1)/2},\qquad z>0.
\end{equation*}
For $\mathcal B_\nu=-\partial_s^2-(2\nu+1)s^{-1}\partial_s$, the
Bessel equation and two applications of Green's formula give
\begin{equation*}
 z^4h_\nu v_N(z)=G_N(z),\qquad
 G_N=h_\nu(\mathcal B_\nu^2v_N).
\end{equation*}
There are no boundary terms because $v_N$ has compact
support in the annulus. The function $\mathcal B_\nu^2v_N$
has the same support and its $L^1(d\mu_\nu)$ norm is
bounded by $C_\nu A_4(a)$. Applying the preceding Bessel
bound inside its Hankel integral gives
\begin{equation*}
 |G_N(z)|+|G_N'(z)|
 \lesssim_\nu A_4(a)z^{-(d-1)/2},\qquad z\ge1.
\end{equation*}
Here differentiation only acts on the Bessel factor,
producing the bounded annular factor $s$; it requires
no additional derivative of $v_N$. Since
\begin{equation*}
 (h_\nu v_N)'(z)
 =z^{-4}G_N'(z)-4z^{-5}G_N(z),
\end{equation*}
both $h_\nu v_N(z)$ and its derivative are bounded by
$C_\nu A_4(a)z^{-L}$ for $z\ge1$, where
$L=4+(d-1)/2$. For $0<z\le1$, their defining
integrals give the bound $C_\nu A_4(a)$.
Since $k_N(z)=N^d(h_\nu v_N)(Nz)$, this gives, uniformly in $N$,
\begin{equation*}
 |k_N(z)|+N^{-1}|k_N'(z)|
 \lesssim_\nu A_4(a)N^d(1+Nz)^{-L}.
\end{equation*}
For $\ell=0,1$, the dyadic sum of
$N^{d+\ell}(1+Nz)^{-L}$ is bounded by $C_\nu z^{-d-\ell}$:
the terms $Nz\le1$ form a geometric sum with positive
exponent $d+\ell$, and the terms $Nz>1$ form a convergent
geometric sum because
\begin{equation*}
 L-d-1=\frac{5-d}{2}>0,\qquad 1<d\le4.
\end{equation*}
Consequently, four symbol derivatives give both bounds
\begin{equation*}
 \sum_{N\in2^{\mathbb Z}}|k_N^{(\ell)}(z)|
 \lesssim_\nu A_4(a)z^{-d-\ell},
 \qquad \ell=0,1,\quad z>0.
\end{equation*}

Writing out the two Hankel transforms, the kernel of $h_\nu a_Nh_\nu$
relative to $d\mu_\nu$ is
\begin{equation*}
 K_N(x,y)=\int_0^\infty a_N(\xi)
 \widetilde J_\nu(x\xi)\widetilde J_\nu(y\xi)\,d\mu_\nu(\xi).
\end{equation*}
For $\nu>-1/2$, the Bessel product formula, in our normalization, reads
\begin{equation*}
 \begin{aligned}
 \widetilde J_\nu(x\xi)\widetilde J_\nu(y\xi)
 &=c_\nu\int_0^\pi\widetilde J_\nu(\xi z_\theta)
                    (\sin\theta)^{d-2}\,d\theta,\\
 z_\theta&=(x^2+y^2-2xy\cos\theta)^{1/2},
 \qquad c_\nu=\bigl(2^\nu\sqrt\pi\,\Gamma(\nu+1/2)\bigr)^{-1}.
 \end{aligned}
\end{equation*}
See \cite[Section~11.41, (16), p.~367]{Watson},
with the cylinder function there chosen to be $J_\nu$.
Substituting this identity and using $k_N=h_\nu a_N$ gives
\begin{equation*}
 K_N(x,y)=c_\nu\int_0^\pi k_N(z_\theta)(\sin\theta)^{d-2}\,d\theta.
\end{equation*}
Fubini's theorem applies because $a_N$ has compact support,
$\widetilde J_\nu$ is bounded, and $d-2=2\nu>-1$ makes the angular weight
integrable.
For $x\ne y$ and $\ell=0,1$,
\begin{equation*}
 \int_0^\pi z_\theta^{-d-\ell}(\sin\theta)^{d-2}\,d\theta
 \lesssim_\nu
 \frac1{|x-y|^{1+\ell}(x+y)^{d-1}}.
\end{equation*}
Indeed, when $|x-y|\ge(x+y)/2$ use $z_\theta\ge|x-y|$.
Otherwise $x\sim y$ and
$z_\theta^2\sim |x-y|^2+(x+y)^2\theta^2$ on $(0,\pi/2)$.
The substitution $\theta=|x-y|u/(x+y)$ reduces this part to an
integral bounded at zero by $u^{d-2}$ and at infinity by $u^{-2-\ell}$.
The other angular half is bounded directly.

The preceding dyadic sums, the angular integral bound and
$|\partial_xz_\theta|,|\partial_yz_\theta|\le1$ yield an
absolutely convergent kernel sum $K=\sum_NK_N$ away from
the diagonal, with
\begin{equation*}
 \begin{aligned}
 |K(x,y)|&\lesssim_\nu
       \frac{A_4(a)}{|x-y|(x+y)^{d-1}},\\
 |\partial_xK(x,y)|+|\partial_yK(x,y)|
       &\lesssim_\nu
       \frac{A_4(a)}{|x-y|^2(x+y)^{d-1}}.
 \end{aligned}
\end{equation*}
We apply Calder\'on--Zygmund theory on the metric measure
space $((0,\infty),|\cdot|,\mu_\nu)$. Its balls are
$B(x,\rho)=(\max\{0,x-\rho\},x+\rho)$, and direct
integration gives
\begin{equation*}
 \mu_\nu(B(x,\rho))
 =\frac{(x+\rho)^d-\max\{0,x-\rho\}^d}{d}
 \simeq_\nu\rho(x+\rho)^{d-1}.
\end{equation*}
In particular, the measure is doubling:
$\mu_\nu(B(x,2\rho))\lesssim_\nu\mu_\nu(B(x,\rho))$.
Put $\delta=|x-y|$. Since $x+\delta\simeq x+y$,
the preceding kernel bounds can be rewritten as
\begin{equation*}
 \begin{aligned}
 |K(x,y)|&\lesssim_\nu
 \frac{A_4(a)}{\mu_\nu(B(x,\delta))},\\
 |\partial_xK(x,y)|+|\partial_yK(x,y)|
 &\lesssim_\nu
 \frac{A_4(a)}{\delta\,\mu_\nu(B(x,\delta))}.
 \end{aligned}
\end{equation*}
These are the kernel size and gradient conditions used
in Calder\'on--Zygmund theory on this space; see
\cite[Section~2]{BCN}, with their parameter
$\lambda=\nu+1/2$. For example, if $|x-x'|<\delta/2$,
the mean value theorem gives the smoothness estimate
\begin{equation*}
 |K(x,y)-K(x',y)|\lesssim_\nu
 A_4(a)\frac{|x-x'|}{\delta}
 \frac{1}{\mu_\nu(B(x,\delta))}.
\end{equation*}
The analogous estimate holds in the second variable.
The operator $h_\nu a h_\nu$ has $L^2(d\mu_\nu)$ norm at most
$\|a\|_\infty$. Its dyadic partial sums converge in $L^2$ by the
spectral representation; for separated compact supports the kernel
bounds also give dominated convergence. Thus $K$ is its off-diagonal
kernel. The weighted Calder\'on--Zygmund theorem gives boundedness on
$L^r(w\,d\mu_\nu)$ for $w\in A_r(d\mu_\nu)$, with norm controlled by
the displayed kernel constants and the
$L^2(d\mu_\nu)$ operator norm of $h_\nu a h_\nu$; see
\cite[Section~2, property~(M1)]{BCN}.
For $w(x)=x^\gamma$, the $A_r$ condition is exactly
$-d<\gamma<d(r-1)$, proving the claim \eqref{eq:claim}.

To return to the original normalization, let $U_\nu g(x)=x^\nu g(x)$.
Direct substitution gives
\begin{equation*}
 \begin{aligned}
 \mathcal H_\nu a\mathcal H_\nu&=U_\nu h_\nu a h_\nu U_\nu^{-1},\\
 \|U_\nu g\|_{L^r(x\,dx)}
     &=\|g\|_{L^r(x^{\nu(r-2)}d\mu_\nu)}.
 \end{aligned}
\end{equation*}
The power-weight conditions reduce to $\nu r>-2$ and $(\nu+2)r>2$.
For $\nu=-\alpha$ these give precisely
$2/(2-\alpha)<r<2/\alpha$; for $\nu=1-\alpha>0$ they hold for every
$1<r<\infty$. Consequently, for both component orders in the stated range,
\begin{equation*}
 \|\mathcal H_\nu a\mathcal H_\nu f\|_{L^r(x\,dx)}
 \lesssim_{\alpha,r} A_4(a)\|f\|_{L^r(x\,dx)}.
\end{equation*}

For a finite set $F\subset\mathbb Z$ and signs $\epsilon_j=\pm1$, take
\begin{equation*}
 a(\xi)=\sum_{j\in F}\epsilon_j\varphi(2^{-j}\omega(\xi)).
\end{equation*}
Its seminorm $A_4(a)$ is bounded independently of $F$
and the signs.
Indeed, for $D=\xi\partial_\xi$ and $b=\xi^2/(m^2+\xi^2)\in[0,1]$,
one has $D=b\,\omega\partial_\omega$ on functions of $\omega$ and
$Db=2b(1-b)$. Repeated differentiation gives bounded polynomials in
$b$ multiplying Euler derivatives of $\varphi$, and only a fixed number
of indices $j$ occur at each $\xi>0$.

Fix either order, set
$T_j=\mathcal H_\nu\varphi(2^{-j}\omega)\mathcal H_\nu$, and use
$\|\cdot\|_r$ for the norm in $L^r(x\,dx)$. Khintchine's inequality gives
\begin{equation*}
  \textstyle
 \left\|\left(\sum_{j\in F}|T_jf|^2\right)^{1/2}\right\|_r
 \lesssim_r
 \left(\mathbb E\left\|\sum_{j\in F}\epsilon_jT_jf\right\|_r^r
 \right)^{1/r} \lesssim_{\alpha,r,\varphi}\|f\|_r.
\end{equation*}
Letting $F$ increase to $\mathbb Z$ proves the upper bound for
$S_\varphi f=(\sum_j|T_jf|^2)^{1/2}$.
Choose a real smooth annular cutoff $\psi=1$ on
$\operatorname{supp}\varphi$, and define $\widetilde T_j,S_\psi$
in the same way. The same upper bound holds for $S_\psi$ at $r'$,
since the stated range is invariant under conjugation.
For $f\in L^2\cap L^r$ and $g\in L^2\cap L^{r'}$, the partition
identity and self-adjointness give
\begin{equation*}
  \textstyle
 |\langle f,g\rangle|
 =\left|\sum_j\langle T_jf,\widetilde T_jg\rangle\right| 
 \le\int_0^\infty S_\varphi f(x)S_\psi g(x)x\,dx
 \lesssim\|S_\varphi f\|_r\|g\|_{r'}.
\end{equation*}
Duality and density prove the lower bound.
The two scalar estimates yield the proposition by
\eqref{eq:massive-functional-calculus-zero-mode}, since the angular
factors have modulus one.
\end{proof}

For the remaining arguments, fix $\varphi$ as in
Proposition~\ref{prop:square function-singular-massive}, and write
$\varphi_j(\lambda)=\varphi(2^{-j}\lambda)$.
Choose $\widetilde\varphi\in C_c^\infty((1/4,4))$ equal to one on
$\operatorname{supp}\varphi$, and put
$\widetilde\varphi_j(\lambda)=\widetilde\varphi(2^{-j}\lambda)$.
We first record the fixed-shell estimate used both here and in the
weighted endpoint argument.

\begin{lemma}[Stability of the fixed-shell estimates]
\label{lem:fixed-shell-stability}
Let $I\Subset(0,\infty)$ be a compact interval and let $\omega_\tau$ be
real smooth phases with $0<c\le|\omega_\tau'|\le C$ on $I$ and uniformly
bounded derivatives. Let $a_\tau$ be smooth amplitudes supported in $I$,
with uniform derivative bounds. For $\nu\in[-1/2,1]$, put

\begin{equation*}
\mathcal E_{\nu,\tau}g(t,r)
=\int_I e^{it\omega_\tau(\lambda)}J_\nu(r\lambda)
a_\tau(\lambda)g(\lambda)\lambda\,d\lambda.
\end{equation*}

For $2\le p,q\le\infty$ and dyadic $R>0$,

\begin{equation}\label{eq:fixed-shell-spacetime}
\|\mathcal E_{\nu,\tau}g\|_{L_t^pL^q((R,2R),rdr)}
\lesssim\|g\|_{L^2(I,\lambda d\lambda)}
\begin{cases}
R^{\nu+2/q},&R\le1,\\
R^{1/q+1/\min\{p,q\}-1/2},&R\ge1.
\end{cases}
\end{equation}

The constants are uniform in $\tau$ and $R$; as usual, $1/\infty=0$.

\end{lemma}

\begin{proof}
The change of variables $\sigma=\omega_\tau(\lambda)$ has a uniformly
bounded Jacobian and inverse Jacobian. Hence Hausdorff--Young and
H\"older on the fixed interval give, for any bounded multiplier $b$,
\begin{equation*}
\left\|\int_I e^{it\omega_\tau(\lambda)}b(\lambda)g(\lambda)
\lambda\,d\lambda\right\|_{L_t^p}
\lesssim\|b\|_{L^\infty(I)}\|g\|_{L^2(I,\lambda d\lambda)}.
\end{equation*}

For $R\le1$, insert the convergent power series for $J_\nu(r\lambda)$.
Apply this time estimate term by term and then sum the absolutely
convergent series. The radial $L^q(rdr)$ norm of each term on $(R,2R)$
has the required leading factor $R^{\nu+2/q}$.

For $R\ge1$, the Bessel and derivative bounds give

\begin{equation*}
\|\partial_r^\ell\mathcal E_{\nu,\tau}g(\cdot,r)\|_{L_t^p}
\lesssim r^{-1/2}\|g\|_2,\qquad \ell=0,1.
\end{equation*}

If $p\ge q$, Minkowski's inequality yields $R^{2/q-1/2}$.
If $p\le q$ and $p<\infty$, use
$W^{1,p}(R,2R;dr)\hookrightarrow L^q(R,2R;dr)$, whose constant is
uniform for $R\ge1$, and then Fubini. This gives
$R^{1/q}R^{1/p-1/2}$. The case $p=q=\infty$ follows directly from
the pointwise Bessel bound. These are exactly the two powers in
\eqref{eq:fixed-shell-spacetime}.
\end{proof}

\begin{proposition}[Dyadic spacetime estimate in the singular channel]
\label{prop:dyadic-strichartz-singular-massive}
Let $0<\alpha<1/2$. Assume either $(p,q)=(\infty,2)$ or
$2<p\le\infty$, $2<q<2/\alpha$ and $1/p+1/q<1/2$.
Then, with $s=1-1/p-2/q$,
\begin{equation*}
 \|\varphi_j(|\mathcal D_A|)e^{it\mathcal D_A}P_0f\|_{L_t^pL_x^q}
 \le C\,2^{js}
 \|\widetilde\varphi_j(|\mathcal D_A|)P_0f\|_{[L^2]^2},
 \qquad j\in\mathbb Z.
\end{equation*}
\end{proposition}

\begin{proof}
The energy case follows from the spectral theorem.
First assume $4<q<2/\alpha$.
Use \eqref{eq:single-bessel-flow}, with the normalized coefficients
from \eqref{eq:momentum-unitarity}. Each component is a sum of two
terms with Bessel order $\nu=-\alpha$ or $\nu=1-\alpha$ and a bounded,
uniformly smooth rescaled coefficient. The enlarged cutoff may be
inserted in the data because it is one on the support of $\varphi_j$.
In the stated range,
\begin{equation*}
 \nu+\frac2q>0,\qquad
 \frac1q+\frac1{\min\{p,q\}}-\frac12<0.
\end{equation*}
For the second inequality, use admissibility when $p\le q$ and $q>4$
when $p\ge q$. Thus Minkowski's inequality and the geometric annular
sums in Lemma~\ref{lem:fixed-shell-stability} give a global
fixed-shell spacetime bound.

If $2^j<m/2$, the localized operator vanishes. For $2^j\ge4m$,
set $\xi=2^j\eta$ and $\tau=m2^{-j}\le1/4$.
The phase is $2^jt\sqrt{\eta^2+\tau^2}$, whose first derivative in
$\eta$ is bounded above and below on a fixed compact annulus,
with all higher derivatives uniformly bounded.
Rescaling time and space by $2^j$ in the fixed-shell bound gives
$2^{j(1-1/p-2/q)}=2^{js}$, and unitarity gives the asserted input norm.

For the finitely many remaining energy scales $m/2\le2^j<4m$,
decompose momentum into dyadic shells $\xi\simeq\lambda$,
$\lambda\lesssim_m1$. Remove the constant phase $e^{\pm imt}$ and set
$\xi=\lambda\eta$, $y=\lambda x$ and $T=\lambda^2t$. The rescaled phase is
\begin{equation*}
 \Omega_\lambda(\eta)
 =\frac{\sqrt{m^2+\lambda^2\eta^2}-m}{\lambda^2},
 \qquad
 \Omega_\lambda'(\eta)
 =\frac{\eta}{\sqrt{m^2+\lambda^2\eta^2}}.
\end{equation*}
Its derivative is bounded above and below on the fixed annulus,
and its higher derivatives are uniformly bounded for
$\lambda\lesssim_m1$. The same fixed-shell bound therefore gives
\begin{equation*}
 \|u_\lambda\|_{L_t^pL_x^q}
 \lesssim_m\lambda^\delta\|g_\lambda\|_2,
 \qquad \delta=1-\frac2p-\frac2q>0.
\end{equation*}
Here $\|g_\lambda\|_2$ is the joint $L^2(\xi\,d\xi)$ norm of the
localized normalized coefficients on the two branches.
Finite overlap and Cauchy--Schwarz imply
\begin{equation*}
 \begin{aligned}
 \left\|\sum_\lambda u_\lambda\right\|_{L_t^pL_x^q}
 &\lesssim_m
 \left(\sum_{\lambda\lesssim_m1}\lambda^{2\delta}\right)^{1/2}
 \left(\sum_\lambda\|g_\lambda\|_2^2\right)^{1/2}\\
 &\lesssim_m
 \|\widetilde\varphi_j(|\mathcal D_A|)P_0f\|_{[L^2]^2}.
 \end{aligned}
\end{equation*}
Since $2^j\simeq m$ on these remaining scales, this has the required
form with the factor $2^{js}$ and a constant depending on the fixed mass.

For $2<q\le4$, choose any $4<q_1<2/\alpha$ and set
\begin{equation*}
 \vartheta=\frac{1/2-1/q}{1/2-1/q_1},\qquad
 p_1=\vartheta p
 \quad\text{with }p_1=\infty\text{ if }p=\infty.
\end{equation*}
Then $0<\vartheta<1$ and $1/p_1<1/2-1/q_1$.
Apply the bound just proved and the energy bound to the linear operator
$f\mapsto\varphi_j(|\mathcal D_A|)e^{it\mathcal D_A}P_0f$,
viewed as an operator on $[L^2]^2$.
Interpolation gives the required time and space exponents and the factor
$2^{js}$, because
\begin{equation*}
 \frac1p=\frac{\vartheta}{p_1},\qquad
 \frac1q=\frac{1-\vartheta}{2}+\frac{\vartheta}{q_1},\qquad
 \vartheta\left(1-\frac1{p_1}-\frac2{q_1}\right)
 =1-\frac1p-\frac2q.
\end{equation*}
Finally, replace $f$ by $\widetilde\varphi_j(|\mathcal D_A|)P_0f$:
this leaves the output unchanged and gives the asserted input norm.
\end{proof}

\begin{proof}[Proof of Theorem~\ref{th:strichartz-singular-massive}]
The energy case $(p,q)=(\infty,2)$ follows from the spectral theorem.
By reflection it suffices to assume $0<\alpha<1/2$.
Let $2<q<2/\alpha$, $2<p\le\infty$ and $1/p+1/q<1/2$,
and put $s=1-1/p-2/q$.

For smooth compactly supported spectral data away from the thresholds,
write
\begin{equation*}
 u=e^{it\mathcal D_A}P_0f,\qquad
 u_j=\varphi_j(|\mathcal D_A|)u,\qquad u=\sum_j u_j.
\end{equation*}
Proposition~\ref{prop:square function-singular-massive} and
Minkowski's inequality, followed by
Proposition~\ref{prop:dyadic-strichartz-singular-massive}, give
\begin{equation*}
 \|u\|_{L_t^pL_x^q}
 \lesssim\left(\sum_j\|u_j\|_{L_t^pL_x^q}^2\right)^{1/2} \lesssim
 \left(\sum_j2^{2js}
 \|\widetilde\varphi_j(|\mathcal D_A|)P_0f\|_{[L^2]^2}^2\right)^{1/2}
 \lesssim_m\|P_0f\|_{H_A^s}.
\end{equation*}
The last inequality follows from the $L^2$ spectral calculus and finite
overlap; on the spectrum $|\rho|\ge m>0$, the dyadic weights are
comparable to the inhomogeneous Sobolev weight, with fixed-mass constants.
Density gives the estimate for all data in $H_A^s$.
Reflection proves the other flux range and completes the proof.
\end{proof}

\subsection{Weighted Strichartz estimates for the singular component: proof of
Theorem~\ref{prop:weighted-singular-massive}}
We prove the weighted estimate as in Section~6.2 of \cite{cacdanyinzha}.
First we prove a weighted $L_t^rL_x^\infty$ endpoint; then we interpolate this endpoint
with the unweighted singular Strichartz estimate. The additional feature in the massive
case is the threshold behaviour at $\rho=\pm m$, which produces the additional
condition on the time exponent, as discussed in Remark \ref{rk:threseffect}.

As in the above subsections, we shall restrict the proof to the case
\begin{equation*}
\alpha\in(0,\tfrac12],
\qquad
\gamma=\frac{3\pi}{4}
\end{equation*}
We write $\mathcal D_A=\mathcal D_{A,\gamma}$ and set
\begin{equation*}
        W_\varepsilon(r)=
        \begin{pmatrix}
        (1+r^{-\alpha-\varepsilon})^{-1} & 0\\
        0&1
        \end{pmatrix}.
\end{equation*}
The case $\alpha\in(\frac12,1)$ is obtained by replacing $\alpha$ with
$1-\alpha$ and interchanging the two spinor components. Indeed, for the
distinguished extension in this range the singular branch is built on
$J_{\alpha-1}=J_{-(1-\alpha)}$, and the weight acts on the second component.

\begin{lemma}[Weighted endpoint]
  \label{lem:weighted-endpoint-massive}
  Let $0<\alpha\le\frac12$, $\varepsilon>0$, and let
  \begin{equation*}
          r>r_\alpha:=\frac2{1-\alpha}.
  \end{equation*}
  Then
  \begin{equation*}
    \|W_\varepsilon e^{it\mathcal{D}_A}P_0f\|
    _{L_t^r(\mathbb R;[L_x^\infty(\mathbb R^2)]^2)}
    \le C_{\varepsilon,r}
    \|P_0f\|_{H_A^{1-\frac1r}(\mathbb R^2;\mathbb C^2)}.
  \end{equation*}
  The same estimate also holds for $r=\infty$; in that case the
  right-hand Sobolev exponent is $1$.
\end{lemma}

\begin{proof}
Let us split the analysis of the singular flow into a threshold part and a high frequency part. Let $\chi\in C_c^\infty((0,\infty))$ satisfy
\begin{equation*}
\chi(\lambda)=1 \quad\text{for }m\leq\lambda\leq4m,
\qquad
\operatorname{supp}\chi\subset[m/2,8m].
\end{equation*}

$\bullet$\emph{High frequency sector: $(1-\chi)(|\mathcal{D}_A|)$.}
Use \eqref{eq:single-bessel-flow} and decompose the momentum into shells
$\xi\simeq N$, $N\gtrsim_m1$. Since $|\rho|>4m$ on the support of
$1-\chi$, momentum and energy are comparable. After scaling $\xi=N\lambda$,
the phase is $Nt\sqrt{\lambda^2+(m/N)^2}$ and its first derivative on
the fixed annulus is bounded away from zero uniformly in $N$.
Lemma~\ref{lem:fixed-shell-stability} therefore gives
\begin{equation*}
\|u_{N,\nu}\|_{L_t^rL^\infty(R,2R)}
\lesssim N^{1-1/r}Q_\nu(NR)\|g_N\|_{L^2(\xi d\xi)},
\quad
Q_\nu(x)=
\begin{cases}x^\nu,&x\le1,\\x^{1/r-1/2},&x\ge1,\end{cases}
\end{equation*}
for $\nu=-\alpha$ or $\nu=1-\alpha$. For the positive order, both
dyadic sums of $Q_\nu(NR)$ are uniformly finite because $r>2$.
For the negative order, write $w_\varepsilon(R)=(1+R^{-\alpha-\varepsilon})^{-1}$.
When $NR\le1$,
\begin{equation*}
w_\varepsilon(R)Q_{-\alpha}(NR)
\le R^{\alpha+\varepsilon}(NR)^{-\alpha}
=N^{-\alpha-\varepsilon}(NR)^\varepsilon
\lesssim_m(NR)^\varepsilon.
\end{equation*}
When $NR\ge1$, the same expression is bounded by $(NR)^{1/r-1/2}$.
The resulting matrix, indexed by dyadic $R$ and $N$, has uniformly
summable rows and columns. Using $\ell^2\hookrightarrow\ell^\infty$ in
the radial annuli, Minkowski's inequality for $r\ge2$, and the Schur test,
we obtain
\begin{equation*}
\|W_\varepsilon e^{it\mathcal{D}_A}(1-\chi)(|\mathcal{D}_A|)P_0f\|_{L_t^rL_x^\infty}
\lesssim_m
\left(\sum_N N^{2(1-1/r)}\|g_N\|_{L^2(\xi d\xi)}^2\right)^{1/2}
\lesssim_m\|P_0f\|_{H_A^{1-1/r}}.
\end{equation*}
Here one sums the two spectral branches and components, using
\eqref{eq:momentum-unitarity}. This argument also applies when $r=\infty$.

\medskip

$\bullet$\emph{Threshold sector: $\chi(|\mathcal{D}_A|)$.}
It remains to estimate the threshold contribution. By the spectral representation
of the singular channel, it is enough to consider finitely many terms of the form
\begin{equation}
    \mathcal T h(t,R)
    :=
    w_\varepsilon(R)
    \int_m^{8m}
        e^{\pm it\rho}
        J_{-\alpha}(\xi R)
        h(\rho)\,d\rho,
    \qquad
    \xi=\sqrt{\rho^2-m^2},
    \label{eq:threshold-model}
\end{equation}
where
\begin{equation*}
    w_\varepsilon(R)
    :=
    (1+R^{-\alpha-\varepsilon})^{-1},
\end{equation*}
and $h\in L^2([m,8m],d\rho)$ is one of the threshold spectral
coefficients, possibly multiplied by a bounded function of $\rho$.

Indeed, the only potentially singular contribution is the one containing
$J_{-\alpha}(\xi R)$ with no factor vanishing at $\rho=m$. The other terms
are easier: those containing $J_{1-\alpha}$ involve a Bessel function of
positive order, while a factor $\sqrt{\rho-m}$ gains one power of $\xi$,
since on the threshold region
\begin{equation*}
    \sqrt{\rho-m}\simeq_m \xi,
    \qquad
    \sqrt{\rho+m}\simeq_m 1.
\end{equation*}
It is therefore sufficient to prove
\begin{equation}
    \|\mathcal T h\|_{L_t^rL_R^\infty}
    \lesssim_{m,\alpha,\varepsilon,r}
    \|h\|_{L_\rho^2},
    \qquad
    r>\frac{2}{1-\alpha}.
    \label{eq:threshold-maximal}
\end{equation}
We decompose the momentum variable $\xi$ dyadically near zero. Let
$\eta\in C_c^\infty((1/2,2))$ be such that
\begin{equation*}
    \sum_{N\in 2^{\mathbb Z}}\eta(\xi/N)=1,
    \qquad \xi>0,
\end{equation*}
and set
\begin{equation*}
    h_N(\rho)
    :=
    \eta\left(\frac{\xi}{N}\right)h(\rho),
    \qquad
    \Omega_N
    :=
    \left\{
        \rho\in[m,8m]:
        \frac N2\leq \sqrt{\rho^2-m^2}\leq 2N
    \right\}.
\end{equation*}
Only dyadic numbers $N\lesssim_m1$ occur. Moreover,
\begin{equation}
    \sum_N\|h_N\|_{L_\rho^2}^2
    \lesssim
    \|h\|_{L_\rho^2}^2,
    \label{eq:threshold-orthogonality}
\end{equation}
and, since
\begin{equation*}
    d\rho
    =
    \frac{\xi}{\sqrt{m^2+\xi^2}}\,d\xi,
\end{equation*}
we have
\begin{equation}
    |\Omega_N|\lesssim_m N^2.
    \label{eq:threshold-shell-size}
\end{equation}
Writing
\begin{equation*}
    \mathcal T_N h_N(t,R)
    :=
    w_\varepsilon(R)
    \int_{\Omega_N}
        e^{\pm it\rho}
        J_{-\alpha}(\xi R)h_N(\rho)\,d\rho,
\end{equation*}
we claim that
\begin{equation}
    \|\mathcal T_N h_N\|_{L_t^rL_R^\infty}
    \lesssim_{m,\alpha,\varepsilon,r}
    N^{1-\frac2r-\alpha}\|h_N\|_{L_\rho^2}.
    \label{eq:threshold-dyadic}
\end{equation}
To prove this, put
\begin{equation*}
    y=NR,
    \qquad
    s=\frac{\xi}{N}.
\end{equation*}
On $\Omega_N$ one has $1/2\leq s\leq2$. Define
\begin{equation*}
    U_N(t,y)
    :=
    \mathcal T_Nh_N(t,y/N).
\end{equation*}
Thus
\begin{equation}
    U_N(t,y)
    =
    w_\varepsilon(y/N)
    \int_{\Omega_N}
        e^{\pm it\rho}
        J_{-\alpha}(sy)h_N(\rho)\,d\rho.
    \label{eq:UN-threshold}
\end{equation}
We first consider $0<y\leq y_0$, where $y_0>0$ is fixed sufficiently
small. Using the convergent expansion
\begin{equation*}
    J_{-\alpha}(z)
    =
    z^{-\alpha}
    \sum_{k=0}^\infty c_k z^{2k},
\end{equation*}
we obtain
\begin{equation*}
    U_N(t,y)
    =
    w_\varepsilon(y/N)y^{-\alpha}
    \sum_{k=0}^\infty c_k y^{2k}
    \int_{\Omega_N}
        e^{\pm it\rho}
        s^{2k-\alpha}h_N(\rho)\,d\rho.
\end{equation*}
Furthermore,
\begin{equation}
    \sup_{y>0}
    w_\varepsilon(y/N)y^{-\alpha}
    \lesssim
    N^{-\alpha}.
    \label{eq:threshold-weight}
\end{equation}
Indeed, setting $x=y/N$,
\begin{equation*}
    w_\varepsilon(y/N)y^{-\alpha}
    =
    N^{-\alpha}
    (1+x^{-\alpha-\varepsilon})^{-1}x^{-\alpha},
\end{equation*}
and the last factor is uniformly bounded for $x>0$.

Taking the supremum in $y$ first, and then applying the
Hausdorff--Young inequality in the $\rho$ variable, gives
\begin{equation*}\begin{aligned}
    \left\|
        \sup_{0<y\leq y_0}|U_N(t,y)|
    \right\|_{L_t^r}
    &\lesssim
    N^{-\alpha}
    \sum_{k=0}^\infty
    |c_k|y_0^{2k}
    \left\|
        \int_{\Omega_N}
            e^{\pm it\rho}
            s^{2k-\alpha}h_N(\rho)\,d\rho
    \right\|_{L_t^r}
    \\
    &\lesssim
    N^{-\alpha}
    |\Omega_N|^{\frac12-\frac1r}
    \|h_N\|_{L_\rho^2}.
\end{aligned}\end{equation*}
Here the series is harmless since $s\in[1/2,2]$ and $y_0$ is fixed.
Using \eqref{eq:threshold-shell-size}, we conclude
\begin{equation}
    \left\|
        \sup_{0<y\leq y_0}|U_N(t,y)|
    \right\|_{L_t^r}
    \lesssim_m
    N^{1-\frac2r-\alpha}
    \|h_N\|_{L_\rho^2}.
    \label{eq:threshold-small-y}
\end{equation}
We next consider $y\geq y_0$. Set
\begin{equation*}
    A_N(y,\rho)
    :=
    w_\varepsilon(y/N)J_{-\alpha}(sy).
\end{equation*}
Since $r>2$, the one-dimensional Sobolev embedding
$W^{1,r}(y_0,\infty)\hookrightarrow L^\infty(y_0,\infty)$ yields
\begin{equation}
    \|U_N\|_{L_t^rL_y^\infty(y_0,\infty)}
    \lesssim
    \|U_N\|_{L_{t,y}^r}
    +
    \|\partial_yU_N\|_{L_{t,y}^r}.
    \label{eq:threshold-sobolev}
\end{equation}
For $\ell=0,1$, Hausdorff--Young in $\rho$ gives, pointwise in $y$,
\begin{equation*}
    \|\partial_y^\ell U_N(\cdot,y)\|_{L_t^r}
    \lesssim
    \|
        \partial_y^\ell A_N(y,\cdot)h_N
    \|_{L_\rho^{r'}}.
\end{equation*}
Let $a$ be defined by
\begin{equation*}
    \frac1{r'}=\frac1a+\frac12,
    \qquad\text{that is}\qquad
    \frac1a=\frac12-\frac1r.
\end{equation*}
By H\"older's inequality,
\begin{equation*}
    \|
        \partial_y^\ell A_N(y,\cdot)h_N
    \|_{L_\rho^{r'}}
    \leq
    |\Omega_N|^{\frac12-\frac1r}
    \sup_{\rho\in\Omega_N}
    |\partial_y^\ell A_N(y,\rho)|
    \|h_N\|_{L_\rho^2}.
\end{equation*}

The standard Bessel estimates, uniformly for $1/2\leq s\leq2$, imply
\begin{equation*}
    |J_{-\alpha}(sy)|
    +
    |\partial_yJ_{-\alpha}(sy)|
    \lesssim
    \begin{cases}
        1, & y_0\leq y\leq1,\\
        y^{-1/2}, & y\geq1.
    \end{cases}
\end{equation*}
Moreover,
\begin{equation*}
    \left|
        \partial_y w_\varepsilon(y/N)
    \right|
    \lesssim \frac1y.
\end{equation*}
Consequently,
\begin{equation}
    \left\|
        \sup_{\rho\in\Omega_N}
        \left(
            |A_N(y,\rho)|
            +
            |\partial_yA_N(y,\rho)|
        \right)
    \right\|_{L_y^r(y_0,\infty)}
    \lesssim 1,
    \label{eq:threshold-A-bound}
\end{equation}
uniformly in $N$. Therefore, from
\eqref{eq:threshold-sobolev} and Fubini,
\begin{equation*}\begin{aligned}
    \|U_N\|_{L_t^rL_y^\infty(y_0,\infty)}
    &\lesssim
    |\Omega_N|^{\frac12-\frac1r}
    \|h_N\|_{L_\rho^2}
    \\
    &\lesssim_m
    N^{1-\frac2r}
    \|h_N\|_{L_\rho^2}.
\end{aligned}\end{equation*}
Since $N\lesssim_m1$, this implies
\begin{equation}
    \|U_N\|_{L_t^rL_y^\infty(y_0,\infty)}
    \lesssim_m
    N^{1-\frac2r-\alpha}
    \|h_N\|_{L_\rho^2}.
    \label{eq:threshold-large-y}
\end{equation}
Together with \eqref{eq:threshold-small-y}, this proves
\eqref{eq:threshold-dyadic}.

We can now sum over the dyadic threshold shells. Since
\begin{equation*}
    \delta
    :=
    1-\frac2r-\alpha>0
\end{equation*}
precisely when
\begin{equation*}
    r>\frac{2}{1-\alpha},
\end{equation*}
we obtain, using Cauchy--Schwarz and
\eqref{eq:threshold-orthogonality},
\begin{equation*}\begin{aligned}
    \|\mathcal Th\|_{L_t^rL_R^\infty}
    &\leq
    \sum_N
    \|\mathcal T_Nh_N\|_{L_t^rL_R^\infty}
    \lesssim
    \sum_N N^\delta\|h_N\|_{L_\rho^2}
    \\
    &\lesssim
    \left(\sum_NN^{2\delta}\right)^{1/2}
    \left(\sum_N\|h_N\|_{L_\rho^2}^2\right)^{1/2}
    \lesssim
    \|h\|_{L_\rho^2}.
\end{aligned}\end{equation*}
This proves \eqref{eq:threshold-maximal} for finite $r$. For $r=\infty$,
Cauchy--Schwarz on $\Omega_N$ and the pointwise Bessel bound give
$\|\mathcal T_Nh_N\|_{L_{t,R}^\infty}\lesssim_m N^{1-\alpha}\|h_N\|_2$,
so the same summation applies.

As observed above, all remaining threshold terms are better: the terms
involving $J_{1-\alpha}$ have no singular factor $N^{-\alpha}$, while
those containing $\sqrt{\rho-m}$ gain one additional power of $N$.
Summing the finitely many spinorial contributions and using the
unitarity of the singular relativistic Hankel transform, we conclude
that
\begin{equation*}
    \|W_\varepsilon e^{it\mathcal{D}_A}\chi(|\mathcal{D}_A|)P_0f\|_{L_t^rL_x^\infty}
    \lesssim
    \|\chi(|\mathcal{D}_A|)P_0f\|_{L_x^2}.
\end{equation*}
Finally, since $|\rho|\simeq m$ on the support of $\chi(|\mathcal{D}_A|)$,
\begin{equation*}
    \|\chi(|\mathcal{D}_A|)P_0f\|_{L_x^2}
    \lesssim_m
    \|P_0f\|_{H_A^{1-\frac1r}}.\qedhere
\end{equation*}
\end{proof}

We now conclude the proof of the weighted Strichartz estimate.

\begin{proof}[Proof of Theorem~\ref{prop:weighted-singular-massive}]
As explained above, we prove the result for $0<\alpha<\frac12$. Thus
\begin{equation*}
  q(\alpha)=\frac2\alpha,
  \qquad
  r_\alpha=\frac2{1-\alpha}.
\end{equation*}

Fix $\varepsilon>0$ and assume
\begin{equation*}
  q(\alpha)\le q<\infty,
  \qquad
  r_\alpha<p\le\infty,
  \qquad
  \theta>\frac{\alpha-2/q}{\alpha+\varepsilon}.
\end{equation*}
Choose $2<q_0<q(\alpha)$, with $q_0$ sufficiently close to 
$q(\alpha)$, and set
\begin{equation*}
  \lambda:=1-\frac{q_0}{q}.
\end{equation*}
Then
\begin{equation*}
  \frac1q=\frac{1-\lambda}{q_0}.
\end{equation*}
Since
\begin{equation*}
  \alpha\lambda\longrightarrow\alpha-\frac2q
  <\theta(\alpha+\varepsilon)
  \qquad\text{as }q_0\uparrow q(\alpha),
\end{equation*}
we may choose $q_0$ close enough to $q(\alpha)$ so that
\begin{equation*}
  \alpha\lambda<\theta(\alpha+\varepsilon).
\end{equation*}
We then choose an independent weight parameter $\delta>0$ so small that
\begin{equation*}
  \lambda(\alpha+\delta)<\theta(\alpha+\varepsilon).
\end{equation*}
We write $W_\delta$ for the weight $W_\varepsilon$ with $\delta$ in place
of $\varepsilon$.

We now choose $p_0>2$ and $r>r_\alpha$ such that
\begin{equation*}
  \frac1p=\frac{1-\lambda}{p_0}+\frac{\lambda}{r},
  \qquad
  \frac1{p_0}+\frac1{q_0}<\frac12 .
\end{equation*}
For finite $p$, the supremum of the allowed right-hand side is
\begin{equation*}
B(q_0):=(1-\lambda)\left(\frac12-\frac1{q_0}\right)
+\frac{\lambda}{r_\alpha}
=\frac{1-\alpha}{2}-\frac{1-\alpha q_0/2}{q}.
\end{equation*}
As $q_0\uparrow2/\alpha$, $B(q_0)\uparrow1/r_\alpha$.
Since $1/p<1/r_\alpha$, choose $q_0$ closer to $q(\alpha)$, if necessary,
so that $1/p<B(q_0)$. This decreases $\lambda$ and preserves the
preceding weight inequality.
Multiplying both endpoint reciprocal exponents by
$c=(1/p)/B(q_0)\in(0,1)$ gives the required strict choices of $p_0$ and $r$.

If $p=\infty$, the same argument is understood with $p_0=r=\infty$.
Since $q_0<q(\alpha)$, Theorem~\ref{th:strichartz-singular-massive} gives
\begin{equation*}
  \|e^{it\mathcal{D}_A}P_0f\|_{L_t^{p_0}L_x^{q_0}}
  \lesssim
  \|P_0f\|_{H_A^{1-\frac1{p_0}-\frac2{q_0}}}.
\end{equation*}
On the other hand, Lemma~\ref{lem:weighted-endpoint-massive}, applied
with the weight parameter $\delta$, gives
\begin{equation*}
  \|W_\delta e^{it\mathcal{D}_A}P_0f\|_{L_t^rL_x^\infty}
  \lesssim
  \|P_0f\|_{H_A^{1-\frac1r}}.
\end{equation*}
Interpolating these two estimates yields
\begin{equation*}
  \|W_\delta^\lambda e^{it\mathcal{D}_A}P_0f\|_{L_t^pL_x^q}
  \lesssim
  \|P_0f\|_{H_A^s},
\end{equation*}
where
\begin{equation*}
  s
  =
  (1-\lambda)\left(1-\frac1{p_0}-\frac2{q_0}\right)
  +
  \lambda\left(1-\frac1r\right)
  =
  1-\frac1p-\frac2q.
\end{equation*}
Finally, write $w_a(R)=(1+R^{-\alpha-a})^{-1}$ for the nontrivial
entry of $W_a(R)$. For $0<R\le1$, the choice of $\delta$ gives
\begin{equation*}
  w_\varepsilon(R)^\theta
  \le R^{\theta(\alpha+\varepsilon)}
  \le R^{\lambda(\alpha+\delta)}
  \le 2^\lambda w_\delta(R)^\lambda.
\end{equation*}
For $R\ge1$, we have
\begin{equation*}
  w_\varepsilon(R)^\theta\le1
  \le2^\lambda w_\delta(R)^\lambda.
\end{equation*}
Thus, componentwise,
\begin{equation*}
  W_\varepsilon^\theta\le2^\lambda W_\delta^\lambda.
\end{equation*}
Therefore
\begin{equation*}
  \|W_\varepsilon^\theta e^{it\mathcal{D}_A}P_0f\|_{L_t^pL_x^q}
  \lesssim
  \|P_0f\|_{H_A^s}.
\end{equation*}
This proves the theorem for $0<\alpha<\frac12$.

For $\frac12<\alpha<1$, the same proof applies after replacing $\alpha$ by
$1-\alpha$ and interchanging the two spinor components. In that case
\begin{equation*}
  q(\alpha)=\frac2{1-\alpha},
  \qquad
  r_\alpha=\frac2\alpha,
\end{equation*}
and the weight acts on the second component. This concludes the proof.
\end{proof}

\end{document}